\documentclass[a4paper,10pt]{article}

\usepackage[top=4cm, bottom=3cm, left=3.2cm, right=3.2cm]{geometry}

\usepackage{Macros}

\begin{document}

\title{Uncertainty quantification with risk measures in production planning}

\author{Simone G\"ottlich\footnotemark[1], \; Stephan Knapp\footnotemark[1]}

\footnotetext[1]{University of Mannheim, Department of Mathematics, 68131 Mannheim, Germany 1 (goettlich@uni-mannheim.de, stknapp@mail.uni-mannheim.de)}

\date{\today}

\maketitle

\begin{abstract}
This paper is concerned with a simulation study for a stochastic production network model, where the capacities of machines may change randomly.  
We introduce performance measures motivated by risk measures from finance leading to a simulation based optimization framework for the production 
planning. The same measures are used to investigate the scenario when capacities are related to workers that are randomly not available. This corresponds 
to the study of a workforce planning problem in an uncertain environment.  
\end{abstract}


\maketitle

\section{Introduction}
Uncertainty quantification is currently an active research topic including 
a wide field of applications. In this work, we focus on the numerical evaluation of performance measures for a production network model whose dynamics is stochastic in the sense that machine failures or capacity drops can randomly occur. Measuring the performance of this stochastic production is related to the quantification of uncertainty or risk. In the production context, we can directly relate risk to the profit of an enterprise and therefore focus on the comparison of classical evaluation methods with monetary risk measures from finance.

Having tailored performance measures at hand, production planning considerations and questions of control can be studied. For example, production planning and control of continuous production models have been studied in \cite{armbruster2012modeling,la2010control} and production planning of discrete production models with randomness in \cite{ji2016optimal}. Obviously, there is a high interest in the definition of performance measures that support planning decisions 
and also allow for optimization purposes, see \cite{Kouri2016,Kouri2018}.

The stochastic production network model under consideration has been originally introduced in \cite{GoettlichKnapp2018} and is based on the deterministic production network model from \cite{ApiceGoettlichHertyBenedetto}. There have been previous stochastic extensions in \cite{GoettlichKnapp2017,GoettlichMartinSickenberger} of the deterministic model but in \cite{GoettlichKnapp2018} the dependence of the machine failures on the actual workload of the machine has been introduced leading to a more complex dynamics. In contrast to agent-based models, the continuous equations govern the evolution of 
aggregated quantities (such as the density of goods) and are valid in the case of homogenous mass production. 
By using this model, we can predict future outcomes of production and also machine breakdowns or capacity drops. A simulation based optimization approach
is applied to study the distribution (or routing) of goods within the network, or, how the capacity should be chosen in advance for a fixed time period. To do so, we nee performance measures, i.e.\ measures that translate the stochastic outcome into reasonable quantities, to validate decisions.

Considering monetary quantities, appropriate performance measures have been introduced in the literature of finance \cite{FoellmerSchied2016} and are based on so-called risk measures. Famous examples are the expectation, Value at Risk and Average Value at Risk (also called conditional Value at Risk). They are originally introduced in the finance area to quantify financial risk.
Since risk measures are fairly general, they can be also used in other contexts. In work \cite{Capolei2015}, risk measures have been introduced for the optimization of oil production. Therein, the focus is the combination and comparison of risk measures in optimization formulations. 
However, our major goal is to focus on risk measures, their evaluation and impact for the stochastic production network model. Therefore, we  
examine how the distribution rates for an optimal routing should be chosen on the base of the introduced performance measures in a simulation based optimization.
Furthermore, we introduce a setting, where the capacity is determined as a sum of individual capacities, which can be on or off, respectively. This allows for a numerical analysis of the maximal possible capacity, e.g.\ number of workers, with respect to the performance measures. 

This paper is organized as follows: in section \ref{sec:2} we introduce the stochastic production network model and its extensions. Section \ref{sec:3} is devoted to the performance measures. Numerical simulation results are studied in section \ref{sec:4}, where mainly two cases are considered: the optimal routing and the workforce planning.

\section{Stochastic model equations}\label{sec:2}
The core model is a production network model consisting of a coupled system of partial and ordinary differential equations.
We start with the introduction of the deterministic model and present then two possible stochastic extensions. 
\subsection{The deterministic model} \label{subsec2}
We briefly recall the production network model from \cite{ApiceGoettlichHertyBenedetto,GoettlichHertyKlar2005} and its stochastic extension to a load-dependent model from \cite{GoettlichKnapp2018}. To focus on the main ideas, we restrict on the case of a production network consisting of one queue-processor unit first. 
This means, we consider a processor with an unbounded queue in front, which represents the storage.
We assume that the processor is defined on an interval $(a,b) \subset \R$, i.e., with length $L = b-a$, and use $\rho(x,t)$ as the density of goods at $x \in (a,b)$ and time $t \geq 0$. Here, the interpretation of $x$ can be the spatial position of the goods within the processor or the so-called degree of completion. The dynamics of the density is given by the following hyperbolic partial differential equation
\begin{equation}
\partial_t \rho(x,t) + \partial_x \min\{v \rho(x,t),\mu\}=0, \label{eq:rho}
\end{equation}
where $\mu \geq 0$ is the maximal capacity and $v >0$ the production velocity. If we prescribe a processing time $T^{\text{prod}}$, we have the relation $v = \frac{L}{T^{\text{prod}}}$.
The queue is placed in front of the processor and its length $q$ is modeled by the following ordinary differential equation
\begin{equation}
\partial_t q(t) = g_{\text{in}}(t)-g_{\text{out}}(t), \label{eq:q}
\end{equation}
i.e.\ as the balance between inflow into and outflow. Here, $g_{\text{in}}(t)$ is the inflow into the queue, which can be an externally given inflow $G_{\text{in}}(t)$ in the case if there is no predecessor. If the queue has predecessors, the inflow into the queue is the weighted sum of the outflow out of the incoming processors.
The outflow of the queue can be described as follows: if the queue is non-empty, the processor takes its maximal capacity $\mu$, and if the queue is empty, the processor can take the inflow into the queue, which is bounded by the maximal capacity $\mu$. Summarizing, this reads as
\begin{equation*}
g_{\text{out}}(t) = 
\begin{cases}
\min\{G_{\text{in}}(t),\mu\}, &\text{ if } q(t) = 0,\\
\mu, &\text{ if } q(t)>0.
\end{cases}
\end{equation*}
The coupling of the processor and  the corresponding queue is prescribed by a boundary condition $\rho(a,t) = \frac{g_{\text{out}}(t)}{v}$ and
initial conditions $\rho(x,0) = \rho_0(x) \in L^1((a,b))$, $q(0) = q_0 \in \R_{\geq 0}$ are given.

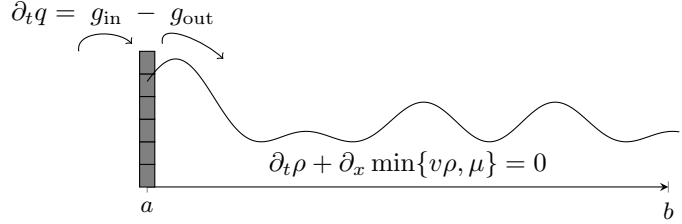
\begin{figure}[H]
\tikzstyle{block} = [rectangle, draw, text width=6em, text centered, sharp corners, minimum height=4em]
\tikzstyle{emptynode} = [draw, ]
\tikzstyle{line} = [draw, -latex']
 \begin{tikzpicture}[font=\small,scale = 1]
\coordinate (origin) at (0,0);
    \begin{axis}[
      xlabel={\normalsize $\partial_t \rho+\partial_x \min\{v\rho,\mu\} = 0$},
      every axis x label/.style={
    at={(ticklabel* cs:0.5)},
    anchor=south,},
    ylabel = {\normalsize $\partial_t q = g_{\text{in}}-g_{\text{out}}$},
          every axis y label/.style={
    at={(ticklabel* cs:1.05)},
    anchor=south,},
      domain =0:1,
      ytick=\empty,
      xtick distance = 1,
      xticklabels = {1,$a$,$b$},
      axis x line=bottom,
      axis y line=none,
    ]
    \end{axis}
  \draw [fill=gray] (-0.1,0) rectangle (0.1,0.3);
   \draw [fill=gray] (-0.1,0.3) rectangle (0.1,0.6);
    \draw [fill=gray] (-0.1,0.6) rectangle (0.1,0.9);
     \draw [fill=gray] (-0.1,0.9) rectangle (0.1,1.2);
     \draw [fill=gray] (-0.1,1.2) rectangle (0.1,1.5);
     \draw [fill=gray] (-0.1,1.5) rectangle (0.1,1.8);
     \node at (-0.45,2.3) (Qeq) {\normalsize $\partial_t q =\; g_{\text{in}}\;-\;g_{\text{out}}$};
     \draw [->] (-0.9,1.8) to [out=100] (-0.2,1.9);
     \draw [->] (0.2,1.9) to [out=100] (1,1.7);
     \draw[domain=0:7.,samples=150,smooth,variable=\x,black] plot ({\x},{0.6+0.4*(sin(\x*100)+1)*(1+cos(\x*120))});
  \end{tikzpicture}
 \caption{Production dynamics on one edge}
 \label{fig:ProductionDynamics}
\end{figure}

Figure \ref{fig:ProductionDynamics} summarizes the modeling equations graphically, where the gray boxes at $a$ represent the queue load following the ODE and the black solid line the density on $(a,b)$ with dynamics described by the hyperbolic PDE.

We now describe the extension to the network case, see \cite{ApiceGoettlichHertyBenedetto}. Suppose that equation \eqref{eq:rho} holds for a density $\rho^e$ on interval $(a^e,b^e)$ with production velocity $v^e>0$ and capacity $\mu^e \geq 0$ for every arc $e \in \cA$ in a directed network $\cG = (\cV,\cA)$, where $\cG$ describes the production network topology with nodes $\cV$ and edges $\cA=\{1,\dots,N\}$. 
We denote by $\delta_v^-$ and $\delta_v^+$ the set of all ingoing and outgoing arcs for every vertex $v \in \cV$. At vertices without any predecessor $v \in V_{\text{in}} = \{v \in \cV \colon \delta_v^-  = \emptyset\}$, we prescribe a time-dependent inflow function $G^v_{\text{in}}(t)$ and for every $v \in \cV$ with $|\delta_v^+|>0$ we assume distribution rates $A^{v,e} \in [0,1], e \in \delta_v^+$, i.e.\ how much of the flow is distributed the subsequent processors. They satisfy $\sum_{e \in \delta_v^+} A^{v,e} (t) = 1$.

\begin{figure}[h]
\centering
\scalebox{1.0}{
\begin{picture}(285,70)(15,-35)
\put(60,0){\vector(1,0){50}}
\put(60,0){\circle*{4}}
\put(80,5){\footnotesize{$1$}}

\put(120,0){\vector(3,2){40}}
\put(120,0){\circle*{4}}
\put(135,18){\footnotesize{$2$}}
\put(122,10){\footnotesize{$\alpha_1$}}

\put(120,0){\vector(3,-2){40}}
\put(120,0){\circle*{4}}
\put(135,-23){\footnotesize{$3$}}
\put(110,-13){\footnotesize{$1-\alpha_1$}}

\put(165,27){\vector(0,-1){50}}
\put(165,27){\circle*{4}}
\put(155,0){\footnotesize{$4$}}

\put(165,27){\vector(3,-2){40}}
\put(165,27){\circle*{4}}
\put(185,18){\footnotesize{$5$}}
\put(174,25){\footnotesize{$\alpha_2$}}

\put(165,-27){\vector(3,2){40}}
\put(165,-27){\circle*{4}}
\put(185,-23){\footnotesize{$6$}}
\put(155,13){\footnotesize{$1-\alpha_2$}}

\put(210,0){\vector(1,0){50}}
\put(210,0){\circle*{4}}
\put(235,5){\footnotesize{$7$}}
\end{picture}
}
\caption{Diamond network with seven processors}
\label{fig:DiamondNetworkTopology} 
\end{figure}
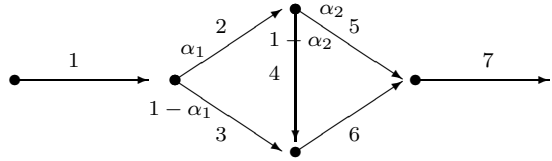

An example of a production network is shown in figure \ref{fig:DiamondNetworkTopology}, where, for example, at node $2$ the distribution rates are given as $A^{2,2} = \ga_1$ and $A^{2,3} = 1-\ga_1$.

We adapt the in- and outflow in and out of the queue of processor $e$ for the network case as follows:
\begin{align*}
g_{\text{in}}^e(t) &= 
\begin{cases}
A^{s(e),e}(t)\sum_{\tilde e \in \delta_{s(e)}^-} \min\{v^{\tilde e} \rho^{\tilde e}(b^{\tilde e},t),\mu^{\tilde e}\} &\text{ if } s(e) \notin V_{\text{in}},\\
G_{\text{in}}^{s(e)}(t) &\text{ if } s(e) \in V_{\text{in}}\text{,}
\end{cases}
\end{align*}
and 
\begin{align*}
g_{\text{out}}^e(t) &= 
\begin{cases}
\min\{g_{\text{in}}^e(t),\mu^e\} &\text{ if } q^e(t) = 0,\\
\mu^e &\text{ if } q^e(t)>0,
\end{cases}
\end{align*}
where $s(e)$ is the starting node of edge $e$.
This deterministic model is well-defined if the initial densities and the inflow functions are of bounded total variation, see \cite{ApiceGoettlichHertyBenedetto}. Moreover, we can define a solution if the inflow and initial densities are in $L^1$ by using an extended solution operator $S^\mu$, i.e.\ $S^\mu_{st}u$ is the solution starting at time $s$ with initial condition $u = (q_0^1,\dots,q_0^N,\rho_0^1,\dots,\rho_0^N)$ at time $t\geq s$ and with capacities $\mu = (\mu^1,\dots,\mu^N)$, see \cite{Bressan2000, ApiceGoettlichHertyBenedetto}. This deterministic flow $S^\mu$ will determine the deterministic evolution between the machine failures or capacity drops.

\subsection{Stochastic extension: load-dependent capacities}

In the following, we introduce the stochastic production network model from \cite{GoettlichKnapp2018}. We assume that the capacities $\mu^e$ can take finitely many non-negative values $\mu^e(i)$ for $i \in \{1,\dots,C^e\}$, $C^e\in \N$ and we introduce the variable $r(t) = (r^1(t),\dots,r^N(t))\in \{1,\dots,C^1\}\times \cdots \times \{1,\dots,C^N\}$ determining the capacities used in the production network at time $t$.
Combining all ingredients leads to deterministic dynamics
\begin{align*}
\Phi_{st} \colon E &\to E,\\
(r_0,q_0,\rho_0) & \mapsto (r(t),q(t),\rho(t)),
\end{align*}
with $$E = \{1,\dots,C^1\}\times \cdots \times \{1,\dots,C^N\} \times \R_{\geq 0}^N \times L^1((a_1,b_1))\times \cdots \times L^1((a_N,b_N)),$$
where
$$r(t) = r_0, \quad (q(t),\rho(t)) = S^{\mu(r_0)}_{st}(q_0,\rho_0).$$
To incorporate random capacity drops in the production network, we follow the theory of piecewise deterministic Markov processes, see e.g. \cite{Davis1984,Jacobsen2006}. Since we have the deterministic evolution between the jump times given by $\Phi$, we only need to specify the intensity $\psi$ at which jumps occur and the distribution of the jumps $\eta$.
To do so, we use rate functions $\gl_{ij}^e$ describing the rate that processor $e$ has a capacity change from $\mu^e(i)$ to $\mu^e(j)$ and assume $$\gl_{ii}^e = \sum_{j\neq i}^{C^e} \gl_{ij}^e.$$
That means, we assume for all $y = (r,q,\rho) \in E$ and $B \in \sigma(E)$
\begin{align*}
\psi(t,y) &= \sum_{e=1}^N \gl^e_{r_er_e}(t,(q_e,\rho_e)),\\
\eta(t,y,B) &= \sum_{e=1}^N\sum_{l\neq r_e}^{C^e} \frac{\gl^e(t,(q^e,\rho^e))}{\psi(t,y)} \varepsilon_{(r_1,\dots,r_{e-1},l,r_{e+1},\dots,r_N,q,\rho)},
\end{align*}
where $\varepsilon_x$ is the Dirac measure with unit mass in $x$ and $\sigma(E)$ is the $\sigma$-algebra on $E$. If the rate functions are continuous with respect to $(t,y)$ and uniformly bounded, then there exists a stochastic process $Y = (Y(t), t\in [0,T])\subset E$ on some probability space \OAP, which is piecewise deterministic between the jumps and follows the deterministic production network equations, see \cite{GoettlichKnapp2018}.

\begin{figure}[htb!]
\begin{picture}(300,100)(20,-10)
\put(50,0){\vector(1,0){250}}
\put(60,-4){\line(0,1){8}}
\put(60,-12){\makebox(5,5){$T_{n-1}$}}
\put(140,-4){\line(0,1){8}}
\put(140,-12){\makebox(5,5){$T_{n}$}}
\put(160,0){\circle{5}}
\put(157,-3.25){\makebox(5,5){X}}
\put(190,0){\circle{5}}
\put(187,-3.25){\makebox(5,5){X}}
\put(205,-4){\line(0,1){8}}
\put(205,-12){\makebox(5,5){$T_{n+1}$}}
\put(300,-12){\makebox(5,5){$t$}}

\put(50,0){\vector(0,1){75}}
\put(50,80){\makebox(5,5){$r^e(t)$}}

\put(60,20){\circle*{5}}
\put(60,26){\makebox(5,5){$Y_{n-1}$}}
\put(60,20){\line(1,0){80}}
\put(140,60){\circle*{5}}
\put(140,66){\makebox(5,5){$Y_{n}$}}
\put(140,60){\line(1,0){65}}
\put(205,40){\circle*{5}}
\put(205,46){\makebox(5,5){$Y_{n+1}$}}
\end{picture}
\caption{Thinning algorithm}
\label{fig:ThinningAlgoPlot}
\end{figure}
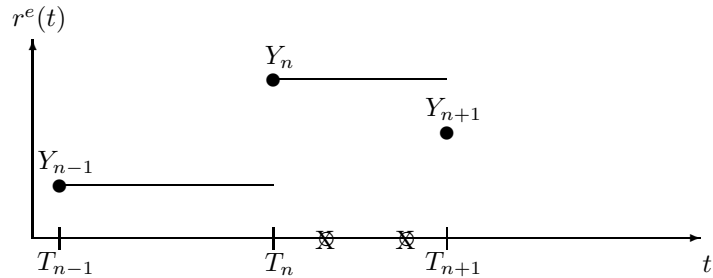

To construct sample paths of the stochastic process $Y$, we use a numerical scheme for the deterministic evolution $\Phi$ and a thinning algorithm for the jump times simultaneously. The approximation of the deterministic evolution consists of a left-sided Upwind scheme for the densities $\rho^e$ and the forward Euler method for the queue-length evolution $q^e$. 
Let $T_n \geq 0$ be the time of the $n$th jump to the value of $Y_n \in E$, then a thinning algorithm produces the next jump time $T_{n+1}$ and post-jump location $Y_{n+1}$ as it is shown schematically in figure \ref{fig:ThinningAlgoPlot}. 

For the description of the procedure we assume a uniform bound $\bar{\lambda}$ on $\psi$. Starting from $T_n$ with the value $Y_n$, we take an exponentially distributed time $\xi_1$ with mean $\bar{\gl}^{-1}$ and use the deterministic evolution to obtain the value $\Phi_{T_n,T_n+\xi_1}Y_n$ at time $T_n+\xi_1$. With an acceptance rejection method, we decide whether a jump is accepted with probability $\tfrac{\Psi(T_n+\xi_1,\Phi_{T_n,T_n+\xi_1}Y_n)}{\bar{\gl}^{-1}}$. If a jump is accepted, the new state of the system $Y_{n+1}$ is produced by the kernel $\eta$. This procedure is repeated until the final time horizon is reached.

\subsection{Stochastic extension: capacities as clusters}\label{subsec:Cluster}

In the following, we study the scenario where the capacity $(\mu^e(t),t\geq 0)$ at processor $e$ consists of $N^e \in \N$ individual capacities. To be more precisely, we assume that $\mu^e(t)$  can be written in the form 
\[\mu^e(t) = \sum_{i=1}^{N^e} X^e_i(t),\] where $(X^e_i(t), t\geq 0)$ is the capacity of cluster part $i$. 
This is important if the production capacity depends on $N^e$ workers that are randomly not available.
We assume that every part of the cluster can be on or off, i.e.\ $X^e_i \in \{0,1\}$, and leads to capacity drops of the capacity process $\mu^e$. In the context of individuals, $X_i^e(t)$ represents whether worker $i$ is available for work or not.

\begin{figure}[H]
\subfigure[Capacities at time $t_1<t_2$]{
 \begin{tikzpicture}[font=\small,scale = 0.6]
    \begin{axis}[
      ybar,
      bar width=20pt,
      xlabel={Processor},
      ylabel={Capacity},
      ymin=0,
      ymax = 10,
      axis x line=bottom,
      axis y line=left,
      enlarge x limits=0.2,
      symbolic x coords={$1$,$2$,$3$,$4$,$5$,$6$,$7$},
      xticklabel style={anchor=base,yshift=-\baselineskip},
      nodes near coords={\pgfmathprintnumber\pgfplotspointmeta}
    ]
      \addplot[fill=white] coordinates {
        ($1$,10)
        ($2$,5)
        ($3$,5)
        ($4$,8)
        ($5$,5)
        ($6$,5)
        ($7$,10)
      };
    \end{axis}
  \end{tikzpicture}
  }
  \subfigure[Capacities at time $t_2>t_1$]{
 \begin{tikzpicture}[font=\small,scale = 0.6]
    \begin{axis}[
      ybar,
      bar width=20pt,
      xlabel={Processor},
      ylabel={Capacity},
      ymin=0,
      ymax = 10,
      axis x line=bottom,
      axis y line=left,
      enlarge x limits=0.2,
      symbolic x coords={$1$,$2$,$3$,$4$,$5$,$6$,$7$},
      xticklabel style={anchor=base,yshift=-\baselineskip},
      nodes near coords={\pgfmathprintnumber\pgfplotspointmeta}
    ]
      \addplot[fill=white] coordinates {
        ($1$,9)
        ($2$,6)
        ($3$,5)
        ($4$,6)
        ($5$,3)
        ($6$,7)
        ($7$,9)
      };
    \end{axis}
  \end{tikzpicture}
  }
  \caption{Worker allocation}
  \label{fig:WorkerAllocation}
\end{figure}
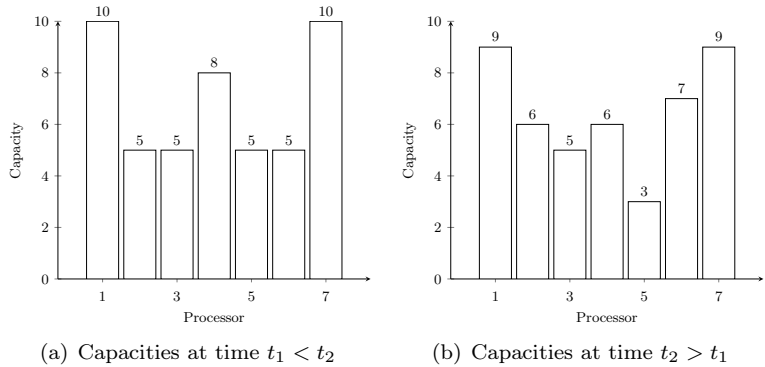

Figure \ref{fig:WorkerAllocation} shows possible realized capacities in the case of the diamond network, cf. figure \ref{fig:DiamondNetworkTopology}, at time $t_1$ and a later time $t_2$. We note that the total capacities are not conserved, i.e.\ workers are either available or not.

The next question is how we can incorporate these ideas into the setting of the model presented in subsection \ref{subsec2}. The mathematical idea is to interpret every $(X_i^e(t), t \geq 0)$ as a continuous time Markov Chain (CTMC) on a common probability space \OAP\ and that they are independent of each other.
The following lemma \ref{lem:sum_CTMC} provides the main tool to numerically evaluate the workforce planning problem
in sub section \ref{subsec:ClusterSIze}. 
To simplify the notation, we neglect the index $e$ in the following.
\begin{lemma}\label{lem:sum_CTMC}
Let,  for some $N\in \N$, a family $(X_1(t),\dots,X_N(t),t\geq 0)$ of independent CTMCs on a common probability space \OAP{} with values in $\{0,1\}$ and $\cQ$-matrix
\begin{align*}
\cQ = 
\begin{pmatrix}
-\gl_0 & \gl_0 \\
\gl_1 & -\gl_1
\end{pmatrix}
\end{align*} 
for $\gl_0, \gl_1>0$ be given.
The stochastic process defined by 
\[X(t)= \sum_{i=1}^N X_i(t)\]
is a CTMC with $\cQ$-matrix satisfying
\begin{align*}
q_{jk} = 
\begin{cases}
-(j\gl_1+(N-j)\gl_0) &\text{ if } k=j\text{,}\\
(N-j)\gl_0	&\text{ if } 		k = j+1\text{,}\\
j\gl_1		&\text{ if } 		k = j-1\text{,}\\
0 		&\text{ else}
\end{cases}
\end{align*}  
for $j,k = 0,\dots,N$. Furthermore, we have
\begin{align}
P(X(t) = j)= \binom{N}{j}P(X_1(t) = 1)^jP(X_1(t) = 0)^{N-j}\text{.} \label{eq:CapDropBino1}
\end{align}
\end{lemma}
The proof can be found in the appendix \ref{app:proof1}.
Note that the proof consists of basically two steps: first proving the Markov property and second using combinatorics to compute the generator $\cQ$ of the process. 

We also have the following remark.
\begin{remark}\label{rem:CapDrop1}
Equation \eqref{eq:CapDropBino1} shows that $X(t)$ is binomially distributed, i.e.\,  \[X(t) \sim \operatorname{Bin}(N,P(X_1(t)= 1)).\] The steady-state distribution is consequently
\begin{align*}
\lim_{t \to \infty} P(X(t) = j) = \binom{N}{j}\left(\frac{\gl_0}{\gl_0+\gl_1}\right)^j\left(\frac{\gl_1}{\gl_0+\gl_1}\right)^{N-j},
\end{align*}
because of \[\lim_{t \to \infty} P(X_1(t) = 1) = \frac{\gl_0}{\gl_0+\gl_1}.\]
\end{remark}
The latter allows for a simple availability analysis of the capacities because the parameters $\gl_0$ and $\gl_1$ are  known from estimations.

Since the entries $q_{ij}$ of the generator $\cQ$ describe the transition rate from state $i$ to $j$ with $i\neq j$, we have 
\begin{align*}
\gl_{ij}^e(t,(q^e,\rho^e)) = 
\begin{cases}
j\gl_1^e+(N^e-j)\gl_0^e &\text{ if } j=i\text{,}\\
(N^e-j)\gl_0^e	&\text{ if } 		j = i+1\text{,}\\
j\gl_1^e		&\text{ if } 		j = i-1\text{,}\\
0 		&\text{ else}
\end{cases}
\end{align*}  
for $i,j \in \{1,\dots,N^e\}$.

This choice embeds this load-independent model into the load-dependent model by the special choice of the rate functions $\gl_{ij}^e(t,(q^e,\rho^e))$ independent of $t,q^e$ and $\rho^e$.

\section{Performance measures}\label{sec:3}
As we have seen, the stochastic production network model introduced allows for random capacities capturing load-dependent failure rates. Since the model is driven by a stochastic process, there is a need for tailored evaluation tools, so-called performance measures, if we have given a set of sample paths. 
This section is devoted to classical performance measures for production models and performance measures based on risk measures motivated from the finance area.

\subsection{Classical performance measures}

Considering sample paths of the production network at every point in time is too detailed in most cases and aggregated quantities are of high interest. According to \cite{GoettlichKnapp2017, GoettlichMartinSickenberger}, the aggregated outflow until time $t\geq 0$ of the complete network can be computed as 
\begin{align*}
G^{\text{net}}_{\text{out}}(t) &= \int_{0}^t \sum_{v \in V_{\text{out}}} \sum_{e \in \delta^-_v}  \min\{v^{ e} \rho^{ e}(b^{ e},s),\mu^{e}(r^e(s))\} ds,
\end{align*}
where $V_{\text{out}} = \{v \in \cV \colon \delta^+_v = \emptyset\}$ are the nodes without a subsequent processor. 
We can also define
\begin{align*}
q^{\text{net}}(t) &= \sum_{e\in \cA} \int_{0}^t q^e(s) ds\text{,}
\end{align*}
as the cumulative sum of all queue-loads up to time $t \geq 0$. Both quantities above are real-valued random variables. Classical performance measures are for example the expectation
$$\E[G^{\text{net}}_{\text{out}}(t)], \quad \E[q^{\text{net}}(t)],$$
and variance 
$$\sigma^2(G^{\text{net}}_{\text{out}}(t)), \quad \sigma^2(q^{\text{net}}(t))$$
of these quantities.

\subsection{Risk measures as performance measures}

If we consider the profit until time $t$ as a random variable $\Pi(t)$ (a functional of $Y$), where $Y$ is a stochastic production network model, 
then we can include monetary aspects.  
We could also use classical measures as the expectation or higher order moments for the random variable $\Pi(t)$ to describe the performance of the production but it turns out that they are not the best choice especially in the context of optimization, see section \ref{sec:4}.
The reason is that a high expected profit can incorporate a high risk that we intend to measure in an appropriate manner. One possibility is the quantification of probability of a bankruptcy, i.e., \[P(\Pi(t)<0).\] One has to keep in mind that this probability does not include any information about the needed surplus to capture ``bad'' events. 
There, we can help us with the so-called Value at Risk and the Average Value at Risk, see e.g.\ \cite{FoellmerSchied2016}, which have been introduced in the context of finance and insurance. They correspond to the class of monetary risk measures, which we introduce in the following definition \ref{def:RiskMeasure} taken from \cite{FoellmerSchied2016}.
\begin{definition}[(Coherent) risk measure]\label{def:RiskMeasure}
Let $\cH$ be a linear space of bounded, real-valued functions containing constants. The mapping $\varrho \colon \cH \to \R$ is called a \emph{monetary risk measure} if
\begin{enumerate}
\item for every $X,\tilde{X} \in \cH$ with $X \leq \tilde{X}$, we have $\varrho(X) \geq \varrho(\tilde{X})$ \hfill (Monotonicity)
\item for every $X \in \cH$, $m \in \R$, it holds that $\rho(X+m) = \rho(X)-m$ \hfill (Cash invariance)
\end{enumerate}
and it is called a \emph{coherent monetary risk measure} if additionally it holds that
\begin{enumerate}
 \setcounter{enumi}{2}
\item for every $\gl \geq 0$, we have $\varrho(\gl X) = \gl \varrho(X)$ \hfill (Positive homogeneity)
\item for every $X,\tilde{X} \in \cH$, it follows that $\varrho(X+\tilde{X}) \leq \varrho(X)+\varrho(\tilde{X})$ \hfill (Subadditivity).
\end{enumerate}
\end{definition}  
Property 1.\ of definition \ref{def:RiskMeasure} states that, if $X$ is interpreted as the profit, the risk of a less profitable company is higher. Assume that we have a fixed surplus of $m$, then the risk is reduced by $m$, see property 2.
Property 3.\ induces a normalization, i.e.\ $\varrho(0) = 0$, and property 4.\ describes risk diversification effect.

One very common risk measure is the Value at Risk ($\operatorname{V@R}_\gl(X)$). It is defined as
\begin{align}
\operatorname{V@R}_\gl(X) = \inf\{m \in \R \colon P(X+m<0)\leq \gl\}
\end{align}
for some level $\gl \in (0,1)$ and a real-valued random variable $X$ on some probability space \OAP; see \cite{FoellmerSchied2016}.
The Value at Risk is simply a quantile and cam be rewritten as
\begin{align}
\operatorname{V@R}_\gl(X) = -q_X^+(\gl). \label{eq:VaR}
\end{align}
with the upper quantile function
\begin{align*}
q^+_X(t) = \sup\{x \in \R \colon  P(X<x) \leq t\}.
\end{align*} 
One can show that $\operatorname{V@R}_\gl$ is a monetary risk measure, which is positive homogeneous but not subadditive. Since the profit of a production network will contain sums of individual costs of machines, we can not guarantee that the Value at Risk incorporates risk diversification effects.
One can easily construct a coherent risk measure from the Value at Risk, which is the Average Value at Risk, and it is defined by 
\begin{align*}
\operatorname{AV@R}_\gl(X) = \frac{1}{\gl}\int_0^\gl \operatorname{V@R}_\gamma(X) d\gamma.
\end{align*}
From the computational point of view, we can estimate the Value at Risk \eqref{eq:VaR} from a sample of profit realizations by an estimation of the quantile function such as \texttt{quantile}\footnotemark[1] in Matlab. This can again be used to compute the Average Value at Risk with, e.g.\  the Matlab function \texttt{integral}\footnotemark[2].

\footnotetext[1]{\href{https://de.mathworks.com/help/stats/quantile.html}{Documentation: https://de.mathworks.com/help/stats/quantile.html}}
\footnotetext[2]{\href{https://de.mathworks.com/help/matlab/ref/integral.html}{Documentation: https://de.mathworks.com/help/matlab/ref/integral.html}}

\section{Computational results}\label{sec:4}

In this section, we numerically investigate the performance measures, their similarities and differences. 
First, we investigate the optimal routing problem and comment on different combinations of distribution rates in a diamond network with respect to the performance measures. Second, we study the impact of cluster sizes, i.e. available workforce, on the performance measures.  

\subsection{Distribution parameter planning in the load-dependent case}\label{subsec:parameterPlanning}

Let the topology of a production network given as a diamond network, see figure \ref{fig:DiamondNetworkTopology}, where $\alpha_1,\alpha_2 \in [0,1]$ are the two distribution parameters, i.e., a percentage of $\alpha_1$ is fed from processor one into queue two ($A^{1,2}(t) = \ga_1$), $1-\alpha_1$ from one to three ($A^{1,3} = 1-\ga_1$) and the same for $\alpha_2$ from processor two to queue five and $1-\ga_2$ to queue four.

As before in section \ref{subsec:ClusterSIze}, we analyze the profit but here for different distribution rates $\ga_1$ and $\ga_2$. To do so, we adopt the profit functional, which is now given as
\begin{align}
\Pi(t) = \int_{0}^t \left(\sum_{v \in V_{\text{out}}} \sum_{e \in \delta^-_v} \min\{v^e\rho^e(b^e,s),\mu^e(r^e(t))\} \cdot p(s) -\sum_{e\in \cA} q^e(s)\cdot C_q^e(s) \right)ds.
\end{align} 
The price of the product is $p(s) \geq 0$ and $C_q^e(s) \geq 0$ is the storage cost at time $s$ for storage $e \in \cA$. 

We assume $p(s) = 1$ and $C_q^e(s) = 0.1$ for every $e = 1,\dots,7$ and $s \in [0,T]$. The queue-processor units all have a production velocity $v^e = 1$ and a length of one, and we start with an empty system at full capacity.
The capacities are given by (ordered by states) $\mu^1 \in \{0,3\}$, $\mu^2,\mu^3 \in \{0,1,2\}$, $\mu^4 \in \{0,1\}$, $\mu^5 \in \{0,2\}$, $\mu^6 \in \{0,1,3\}$, and $\mu^7 \in \{0,2,3\}$.  
To include the load-dependency we follow \cite{GoettlichKnapp2018} and define the Utilization Ratio $\operatorname{UR}$ by
\begin{align*}
\operatorname{UR}^e(r^e,q^e,\rho^e) =  \frac{1}{\max_i\{\mu^e(i)\}(b^e-a^e)} \int_{a^e}^{b^e}\min\{\mu^e(r^e(t)),v^e\rho^e(x,t)\}dx,
\end{align*}
and the Ratio of Work In Progress and the maximal amount of goods in the machine $\operatorname{RWIP}$ as
\begin{align*}
\operatorname{RWIP}^e (r^e,q^e,\rho^e)  = \frac{v^e}{\max_i\{\mu^e(i)\}(b^e-a^e)}\int_{a^e}^{b^e} \rho^e(x,t)dx.
\end{align*}
The rate functions of processors $e \in \{1,4,5\}$ are then given by 
\begin{align*}
\gl_{12}^e(t,q^e,\rho^e) &= \gl^{\text{rep,max,e}}-(\gl^{\text{rep,max,e}}-\gl^{\text{rep,min,e}})\operatorname{RWIP}^e(1,q^e,\rho^e),\\[1ex]
\gl_{21}^e(t,q^e,\rho^e) &= \gl^{\text{down},e}\operatorname{UR}^e(2,q^e,\rho^e),
\end{align*}
with $\gl^{\text{rep,max,e}} = 10$, $\gl^{\text{rep,min,e}} = 4$ and $\gl^{\text{down},e} = 1$.

In the case of three states, i.e., processors 2, 3, 6 and 7, the rate functions read as
\begin{align*}
\gl_{13}^e(t,q^e,\rho^e) &= \gl_{23}^e(t,q^e,\rho^e) \\ &= \gl^{\text{rep,max,e}}-(\gl^{\text{rep,max,e}}-\gl^{\text{rep,min,e}})\operatorname{RWIP}^e(1,q^e,\rho^e),\\[1ex]
\gl_{21}^e(t,q^e,\rho^e) &= \gl^{\text{down},e}\operatorname{UR}(2,q^e,\rho^e),\\[1ex]
\gl_{31}^e(t,q^e,\rho^e) &= \gl_{32}^e(t,q^e,\rho^e)=\gl^{\text{down},e}\operatorname{UR}^e(3,q^e,\rho^e)
\end{align*}
with $\gl^{\text{rep,max,e}} = 10$, $\gl^{\text{rep,min,e}} = 4$ and $\gl^{\text{down},e} = 2$.
This means we consider two different types of processors with two or three capacity states, identical repair rates but different breakdown rates. We can not state an explicit stationary expected capacity as in section \ref{subsec:ClusterSIze} and simulations are needed to determine which combination $(\alpha_1,\alpha_2)$ performs in a ``optimal'' way.

\begin{figure}[h]
\subfigure[Value at Risk]{
\includegraphics[width=0.46\textwidth]{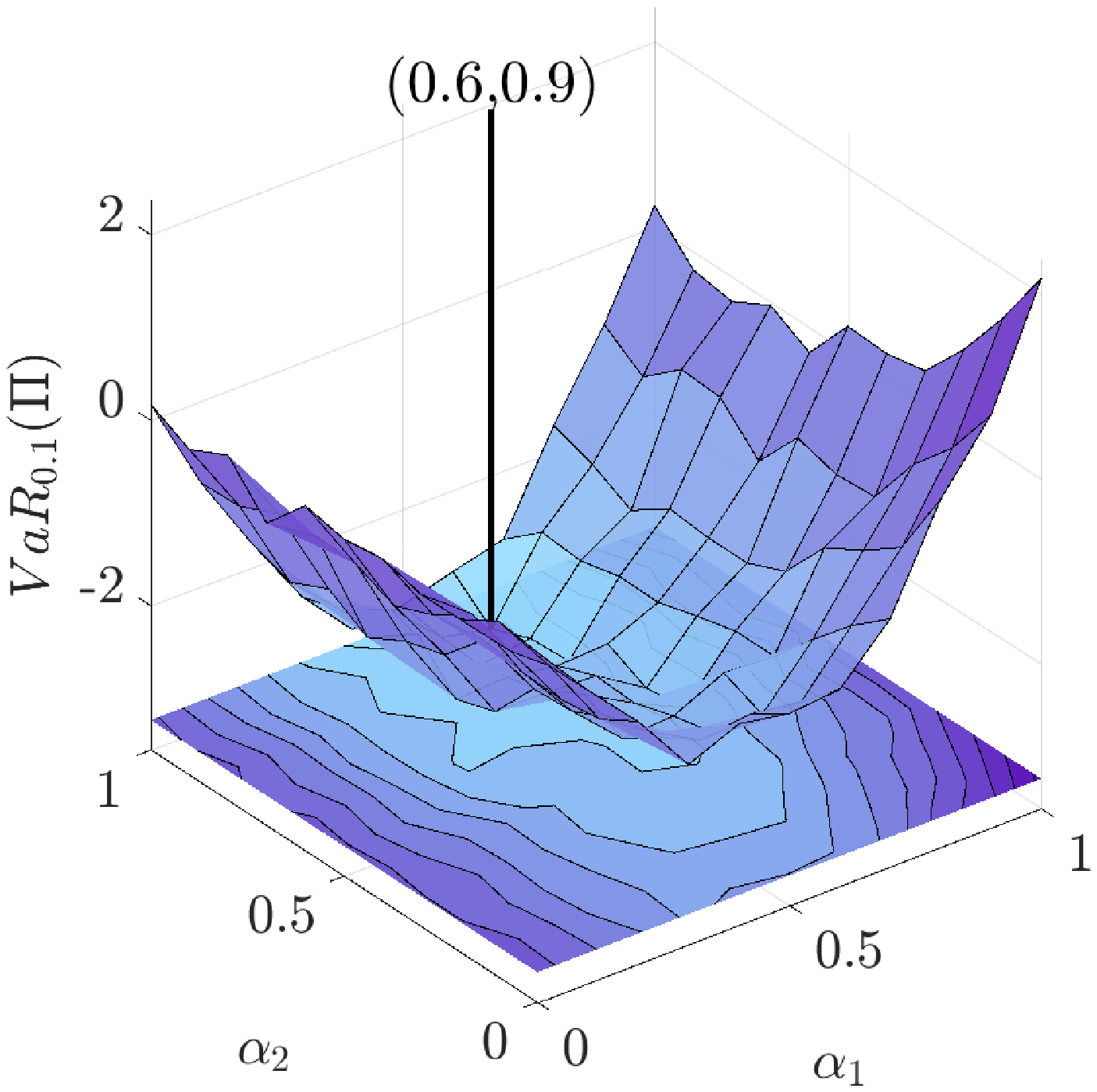}
}
\subfigure[Average Value at Risk]{
\includegraphics[width=0.46\textwidth]{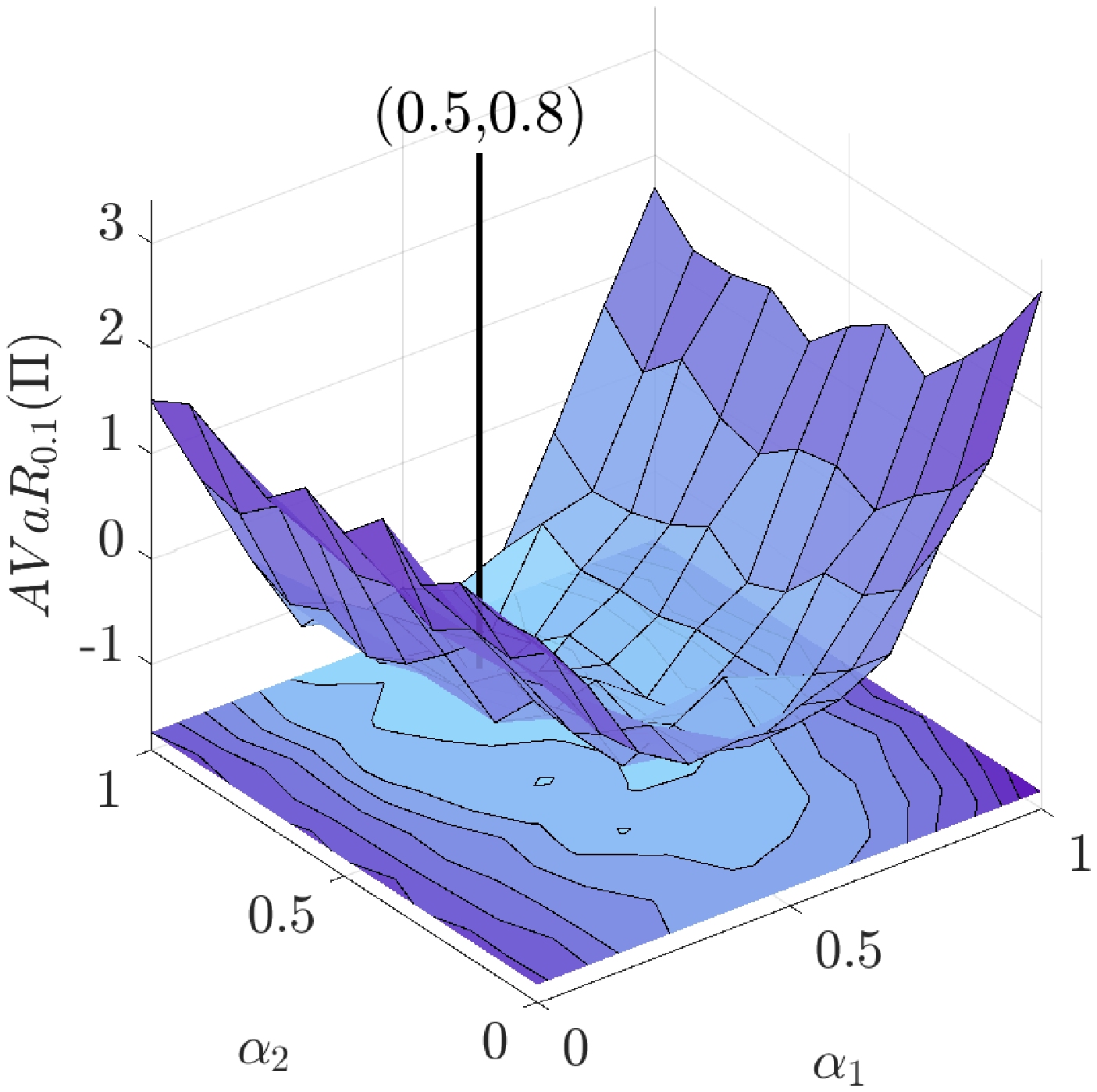}
}
\caption{Value and Average Value at Risk with level $\lambda = 0.1$ for 1000 samples}
\label{fig:LoadDepDiamondProfit2}
\end{figure}

\begin{figure}[h]
\subfigure[Sample mean]{
\includegraphics[width=0.46\textwidth]{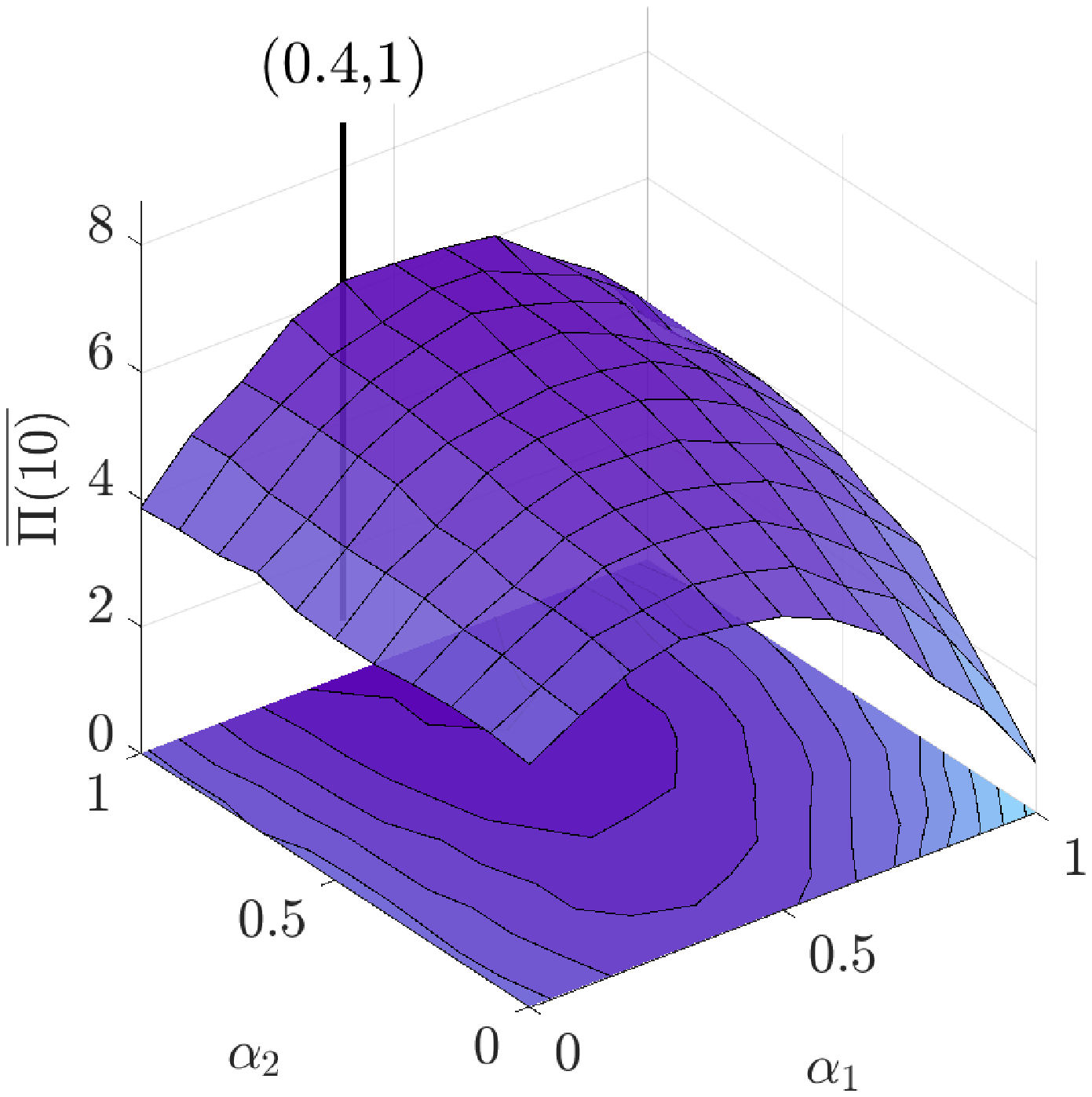}
}
\subfigure[Sample standard deviation]{
\includegraphics[width=0.46\textwidth]{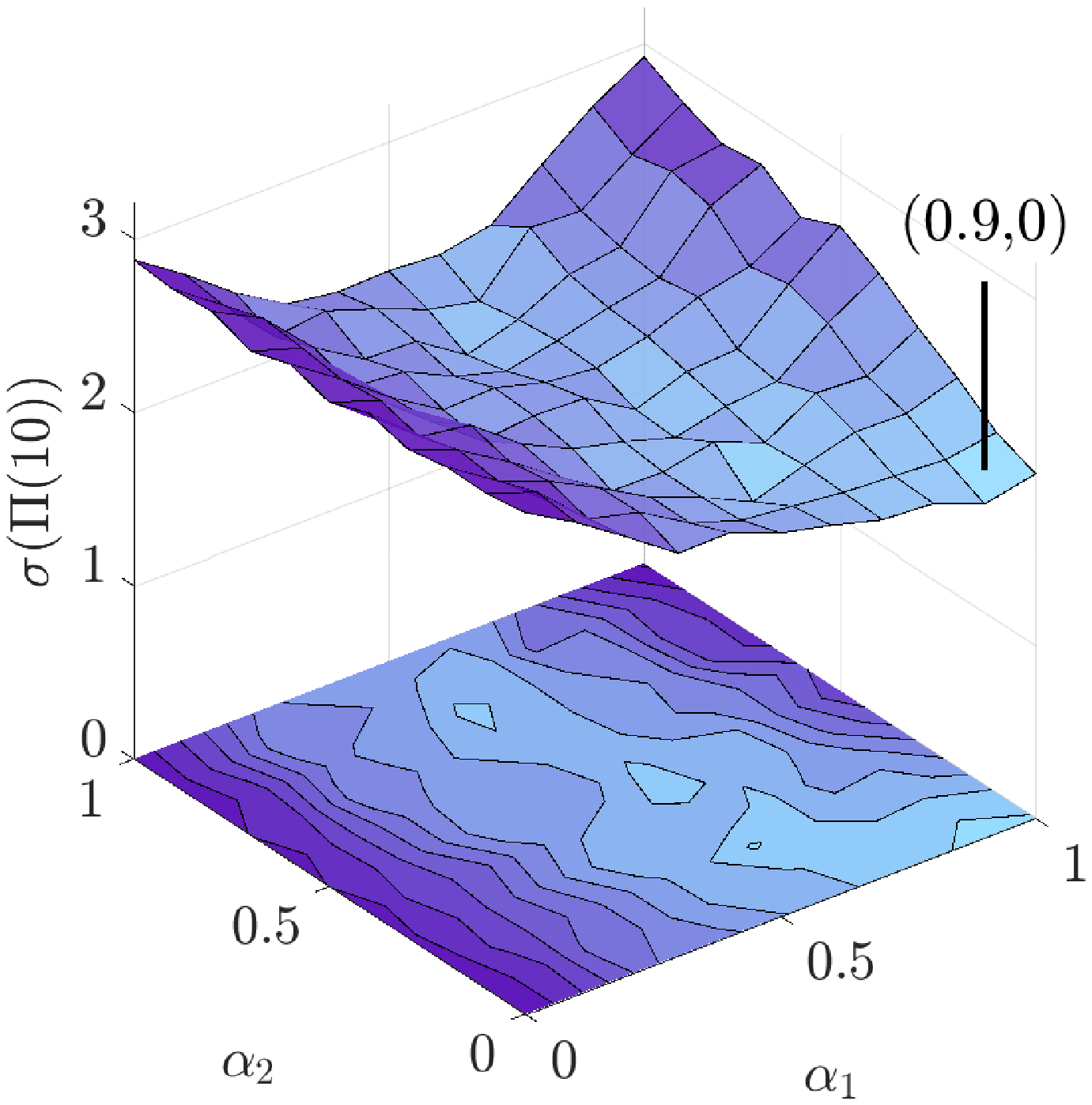}
}\\
\begin{center}
\subfigure[Bankruptcy probability]{
\includegraphics[width=0.46\textwidth]{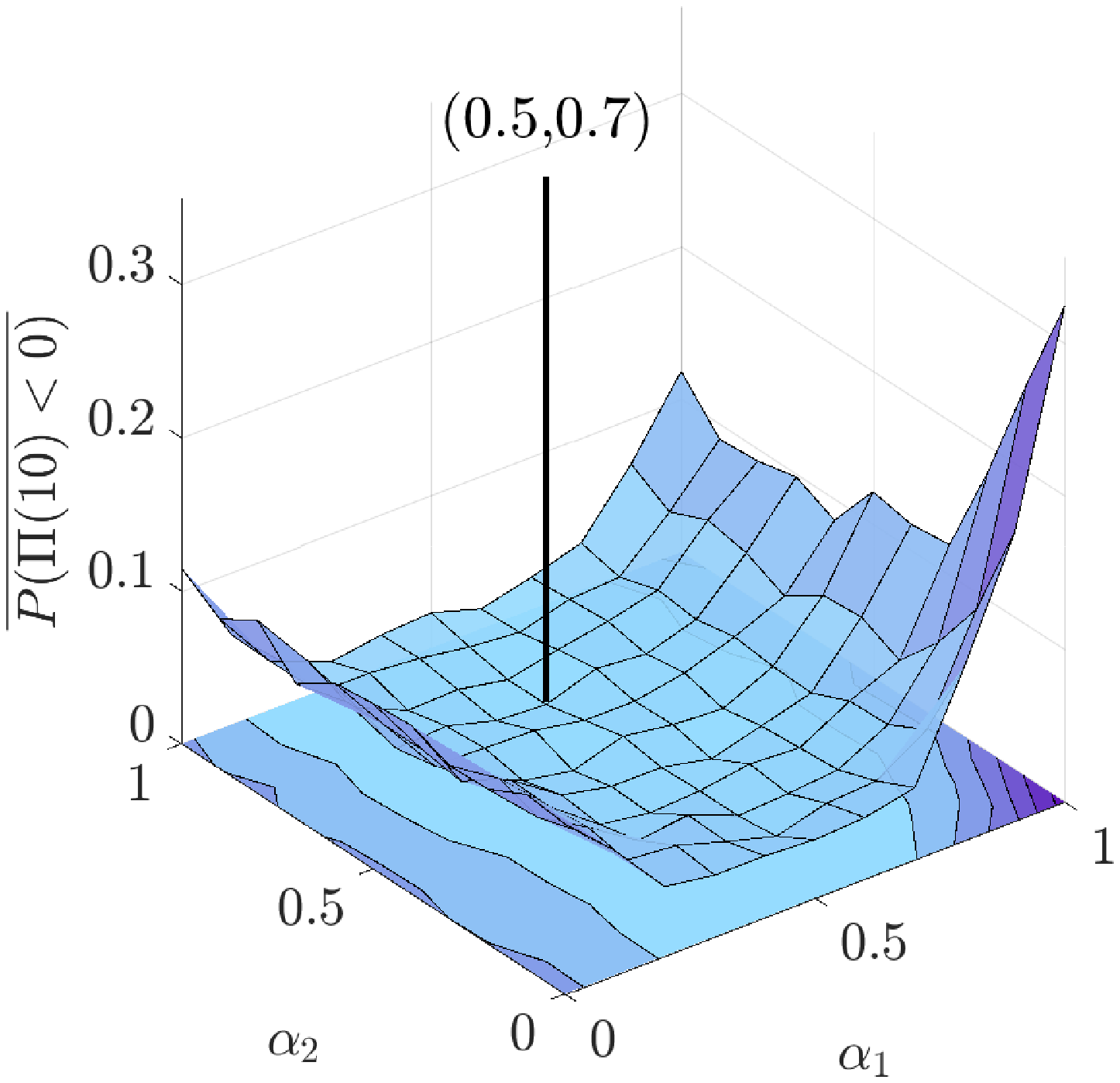}
}
\end{center}
\caption{Sample mean, standard deviation and bankruptcy probability of profit for 1000 samples}
\label{fig:LoadDepDiamondProfit1}
\end{figure}

In figure \ref{fig:LoadDepDiamondProfit1} (a)-(b), the sample mean and standard deviation of the profit depending on different choices of $\alpha_1$ and $\alpha_2$ are shown, where $\ga_1,\ga_2 \in \{0,0.1,\dots,0.9,1\}$ are evaluated. In figures \ref{fig:LoadDepDiamondProfit1} (c) and \ref{fig:LoadDepDiamondProfit2} (a)-(b),  the bankruptcy probability and both the Value at Risk and Average Value at Risk, respectively, are drawn for different combinations of $\ga_1$ and $\ga_2$. The inflow into the network is constant and given by $1.5$, and the time horizon is chosen to be $T = 10$.

We see a strong influence of the distribution parameters on the profit. The vertical black line describes the allocation of $(\alpha_1,\alpha_2)$, which is the best choice for the corresponding evaluation of the profit. In the case of the sample mean, the choice $(\alpha_1,\alpha_2) = (0.4,1)$ leads to the highest value of $6.2087$; see table \ref{tab:LoadDepProfitEval}. In contrast, the sample standard deviation is small for the choice $(\alpha_1,\alpha_2) = (0.9,0)$, with a value of $1.9261$. The lowest bankruptcy probability follows from the choice $(0.5,0.7)$, which is close to the best choices of the Value at Risk $(0.6,0.9)$ and the Average Value at Risk $(0.5,0.8)$ in figure \ref{fig:LoadDepDiamondProfit2}. 
Regarding the different evaluations of the profit sample, the choice of combination $(\ga_1,\ga_2)$ strongly depends on the choice of the evaluation criterion, but there is a tendency to equally distribute at node two and more into processor five than in processor four at node three. Feeding more into processor five can also be motivated by the higher capacity compared to the capacity of processor four. 
The values of the mean capacity, standard deviation, bankruptcy probability and (Average) Value at Risk for all the best choices are listed in table \ref{tab:LoadDepProfitEval}. There, the Average Value at Risk performs well because all other evaluation criteria are not  worse in this case, specifically, because the Average Value at Risk is a coherent risk measure from the theory of risk measures and is thus a reasonable result.

Obviously, the bankruptcy probability is a worse measure in optimizing $(\alpha_1,\alpha_2)$ for this example since the relevant region is very flat and disturbed by the errors from the Monte Carlo simulation. In high contrast, both, the Value at Risk and Average Value at Risk form a nice shaped measure and neglecting the Monte Carlo errors it seems to be convex as well.
\begin{table}[h]
\centering
\small{
\begin{tabular}{l||c|c|c|c|c}
&$(0.4,1)$&$(0.9,0)$&$(0.5,0.7)$&$(0.6,0.9)$&$(0.5,0.8)$\\\hline\hline
$\overline{\Pi(T)}$&$\textbf{6.2087}$&$1.9018$&$5.9841$&$6.1532$&$6.1820$\\\hline
$\sigma(\Pi(T))$&$2.2415$&$\textbf{1.9261}$&$2.1142$&$2.1117$&$2.0602$\\\hline
$\overline{P(\Pi(T)<0)}$&$0.0210$&$0.1890$&$\textbf{0.0120}$&$0.0140$&$\textbf{0.0120}$\\\hline
$\operatorname{V@R}_{0.1}(\Pi(T))$&$-3.0895$&$0.8949$&$-2.9776$&$\textbf{-3.2336}$&$-3.1459$\\ \hline
$\operatorname{AV@R}_{0.1}(\Pi(T))$&$-1.3008$&$1.6910$&$-1.4979$&$-1.5509$&$\textbf{-1.6523}$
\end{tabular}
}
\caption{Comparison of the different profit evaluations}
\label{tab:LoadDepProfitEval}
\end{table}

\subsection{Study of cluster sizes and workforce planning}\label{subsec:ClusterSIze}
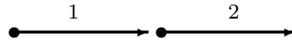
\begin{figure}[H]
\centering
\scalebox{1.0}{
\begin{picture}(185,20)(15,-10)
\put(60,0){\vector(1,0){50}}
\put(60,0){\circle*{4}}
\put(80,5){\footnotesize{$1$}}

\put(115,0){\vector(1,0){50}}
\put(115,0){\circle*{4}}
\put(140,5){\footnotesize{$2$}}
\end{picture}
}
\caption{Serial network with two processors}
\label{fig:Serial2} 
\end{figure}
We analyze the stochastic production network model with two arcs, see figure \ref{fig:Serial2}, where capacities are given as in subsection \ref{subsec:Cluster} and we define the profit $\Pi(t)$ as the cumulative revenue reduced by storage and cluster-size costs until time $t \geq 0$. In formulas, this reads as
\begin{align*}
\Pi(t) = \int_{0}^t \Big(\sum_{v \in V_{\text{out}}} \sum_{e \in \delta^-_v} \min\{v^e\rho^e(b^e,s),\mu^e(r^e(t))\} \cdot p(s)-\sum_{e\in \cA}(q^e(s)\cdot C_q^e(s)+N^e \cdot C_N^e(s))\Big) ds,
\end{align*} 
where we denote with $p(s) \geq 0$ the price of the product, $C_q^e(s) \geq 0$ and $C_N^e(s) \geq 0$ are the storage and cluster size costs, respectively, at time $s$. 
Typical examples for the cluster size cost are maintenance costs or salaries. If $\varrho$ is a monetary risk measure, we can interpret $\varrho(\Pi(T))$ as the risk of the company of the aggregated loss and gains until time $T$. 

In the following, we study an example and consider a production network model in the form of a chain of two processors. 
The capacities of the processors satisfy $\mu^e(t) \in \{0,\dots,N^e\}$ for given cluster sizes $N^1,N^2 \in \N$.
We further assume that the times between failures and the repair times of each part of a cluster are exponentially distributed with mean $\operatorname{MTBF}^1=  80$, $\operatorname{MTBF}^2 = 50$ and $\operatorname{MRT}^1=10$, $\operatorname{MRT}^2 = 20$. This results in the rates
\begin{align*}
\gl_0^1 = \frac{1}{10}, \quad \gl_1^1 = \frac{1}{80}, \quad \gl_0^2 = \frac{1}{20}, \quad \gl_1^2 = \frac{1}{50}.
\end{align*}
We can compute the expected capacities in steady state as
\begin{align*}
\E[\mu^e(\infty)] = \sum_{n=0}^{N^e}n \binom{N^e}{n}\left(\frac{\gl_0^e}{\gl_0^e+\gl_1^e}\right)^n\left(\frac{\gl_1^e}{\gl_0^e+\gl_1^e}\right)^{N^e-n} = N^e \frac{\gl_0^e}{\gl_0^e+\gl_1^e},
\end{align*}
which implies
\begin{align*}
\E[\mu^1(\infty)] = N^1  \frac{8}{9}, \quad \E[\mu^2(\infty)] = N^2  \frac{5}{7},
\end{align*}
and means that the expected capacity of the second processor for $N^1=N^2$ is lower than the capacity of the first one.
We assume a processing velocity $v^e=1$ for both processors, a time step $\Delta t = 1$  and a time horizon of $T = 365$. The length of  each processor is one, and to satisfy the Courant-Friedrichs-Lewy (CFL) stability condition, we choose a spatial step size of $\Delta x=1$.  Given that $G_{\text{in}}^1(t) = 10$ is a constant inflow of ten parts per unit time, we start with an empty production.
The only missing parameters are the cluster sizes $N^1$ and $N^2$ and our goal is to analyze how the profit $\Pi(T)$ changes if we choose different combinations of cluster sizes. 
To evaluate the profit, we assume storage costs $C^1_q = C^2_q = 0.01$, cluster size costs $C_N^1 = 4$, $C_N^2 = 6$ and a price $p = 10.02$.

\begin{figure}[h]
\subfigure[Sample mean]
{
\includegraphics[width=0.46\textwidth]{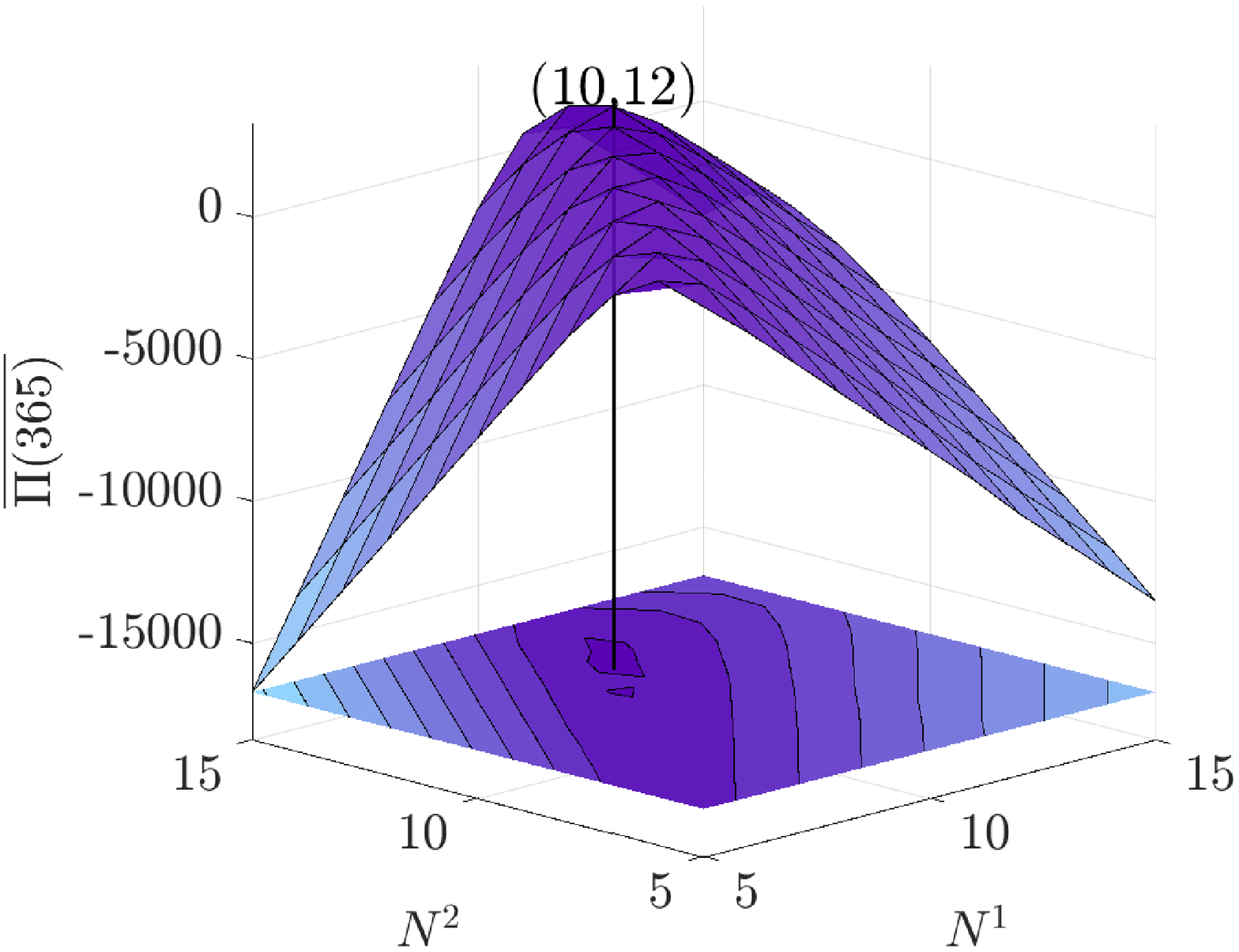}
\label{fig:Optimal_Worker_Mean_PI}
}
\subfigure[Sample standard deviation]
{
\includegraphics[width=0.46\textwidth]{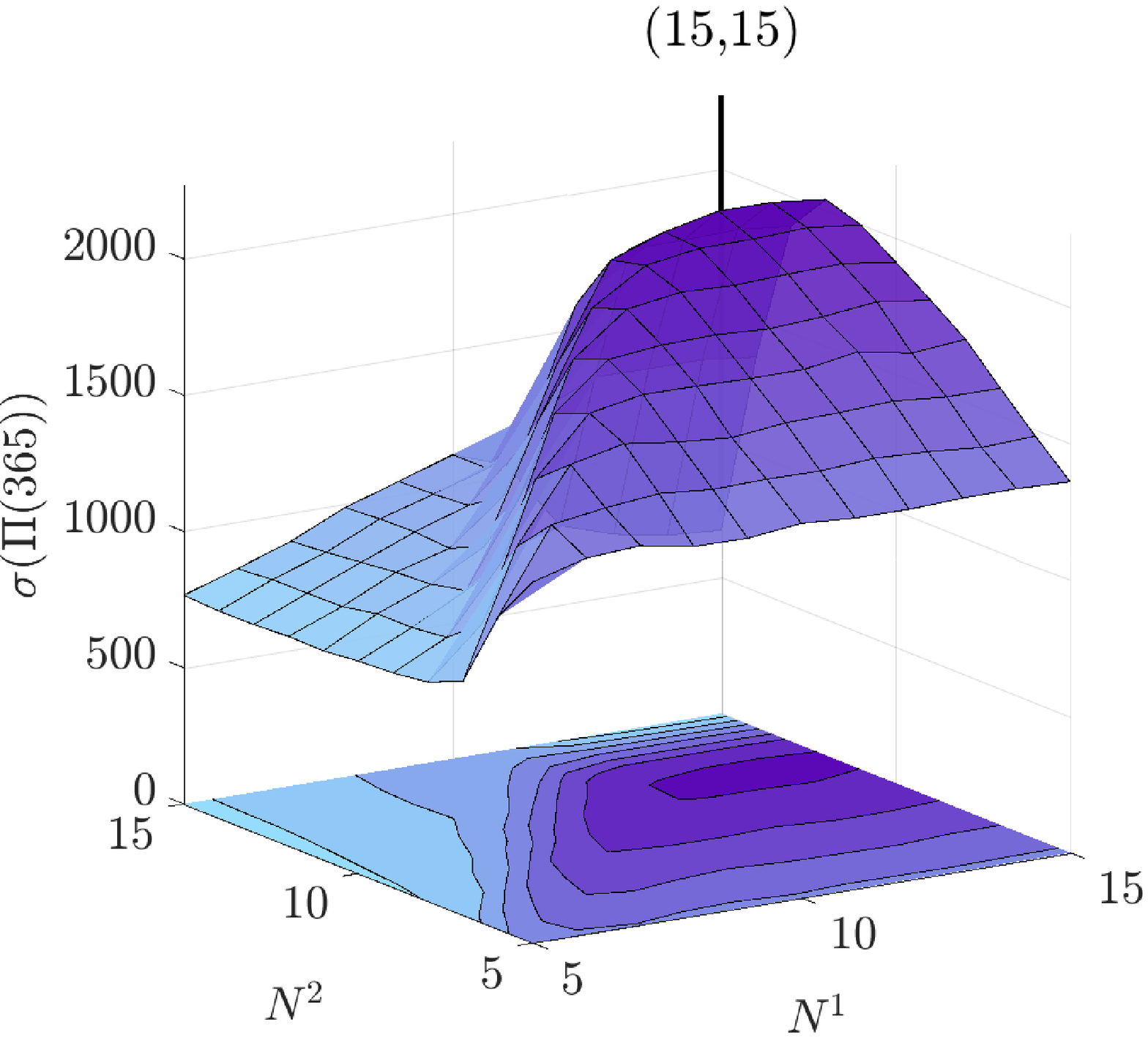}
\label{fig:Optimal_Worker_Std_PI}
}\\
\begin{center}
\subfigure[Estimated bankruptcy probability]
{
\includegraphics[width=0.46\textwidth]{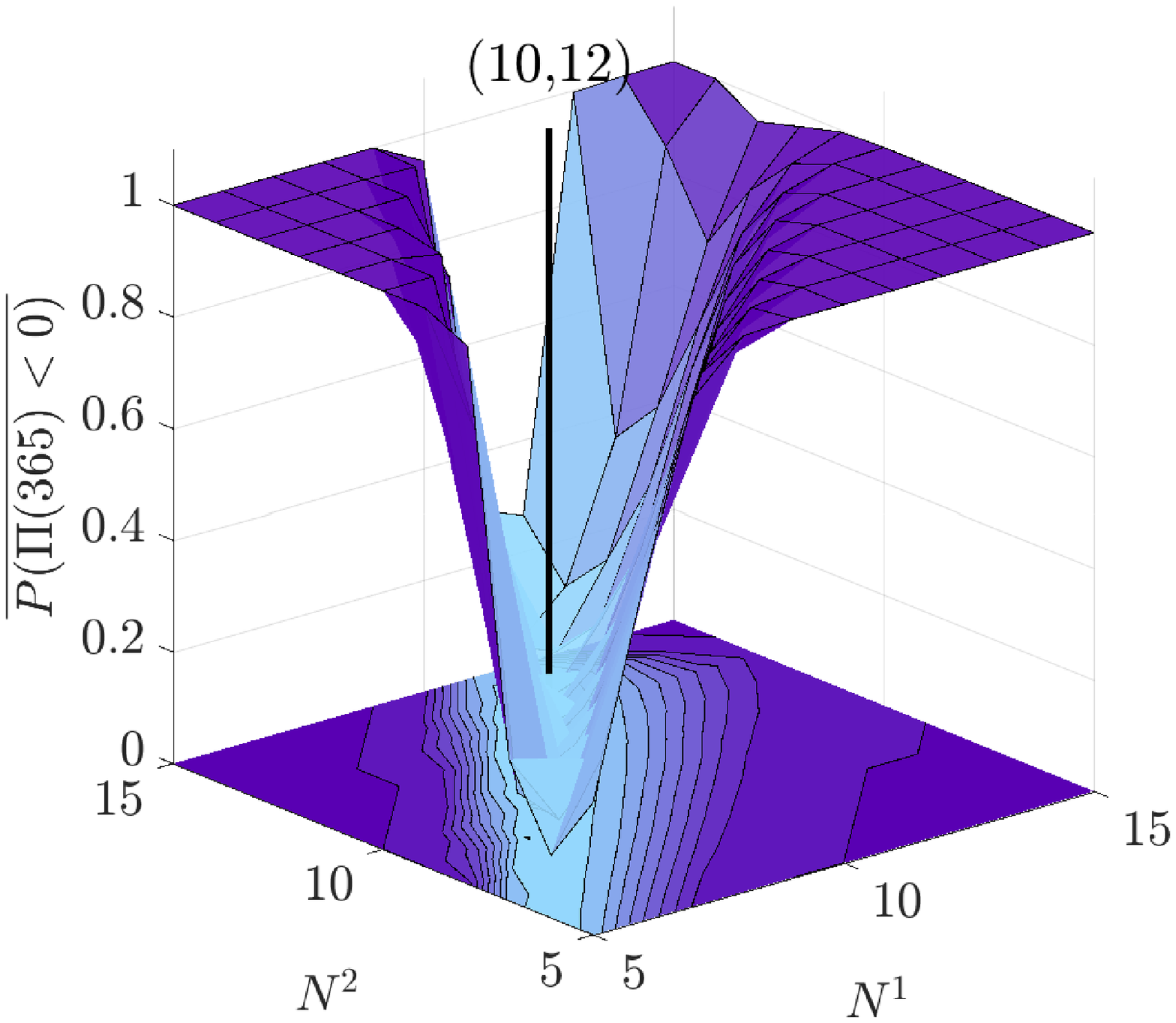}
\label{fig:Optimal_Worker_Bankrupt_Prob}
}
\end{center}
\caption{Profit evaluation for different cluster sizes}
\label{fig:Profit_Classic}
\end{figure}

The following simulation results are based on $M = 10^4$ Monte Carlo samples. Figure \ref{fig:Profit_Classic} contains the sampled values of the expected profit (a), the standard deviation (b) and bankruptcy probability (c). 

In this example, the expected profit strongly depends on the cluster sizes, i.e.\ there are few combinations that lead to a positive expected profit. The highest expected profit is achieved with the choice $N^1 = 10$ and $N^2 = 12$ with a profit of $2.3535 \cdot 10^3$, as table \ref{tab:optim_work} states. 

With respect to the standard deviation, we obtain the combination $N^1 = N^2 = 15$ having the lowest standard deviation of $0.6733\cdot 10^3$; see table \ref{tab:optim_work}. This is reasonable because in this case, the production is less affected by capacity drops due to the large cluster sizes. 

\begin{table}[h]
\centering
\small{
\begin{tabular}{l||c|c|c|c|c}
&$(10,12)$&$(15,15)$&$(10,13)$&$(7,9)$&$(8,10)$\\\hline\hline
$\overline{\Pi}\,[10^3]$&$\textbf{2.3535}$&$-4.0108$&$1.7406$&$1.4955$&$1.8488$\\\hline
$\sigma(\Pi)\,[10^{3}]$&$1.6026$&$\textbf{0.6733}$&$1.2503$&$1.1219$&$1.2804$\\\hline
$\overline{P(\Pi<0)}$&$\textbf{0.0804}$&$1.0000$&$0.0883$&$0.0966$&$0.0831$\\\hline
$\operatorname{V@R}_{0.1}(\Pi)\,[10^3]$&$\textbf{-0.2283}$&$4.8542$&$-0.0995$&$-0.0271$&$-0.1818$\\ \hline
$\operatorname{AV@R}_{0.1}(\Pi)\,[10^3]$&$0.7212$&$5.6145$&$0.6568$&$0.6682$&$\textbf{0.6255}$\\ \hline
$\operatorname{V@R}_{0.01}(\Pi)\,[10^3]$&$1.8642$&$6.5899$&$\textbf{1.5683}$&$1.5774$&$1.6160$\\\hline
$\operatorname{AV@R}_{0.01}(\Pi)\,[10^3]$&$2.6565$&$7.2945$&$2.1790$&$\textbf{2.0632}$&$2.3147$\\
\end{tabular}
}
\caption{Comparison of the different profit evaluations for $(N^1,N^2)$}
\label{tab:optim_work}
\end{table}

The bankruptcy probability in figure \ref{fig:Profit_Classic} (c) shows a tight region in which the probability of going bankrupt is low. The combination with $N^1 = 10$ and $N^2 = 12$ leads to a bankruptcy probability of $0.0804$ and the same combination as the expected profit.

\begin{figure}[h]
\subfigure[Value at Risk]
{
\includegraphics[width=0.46\textwidth]{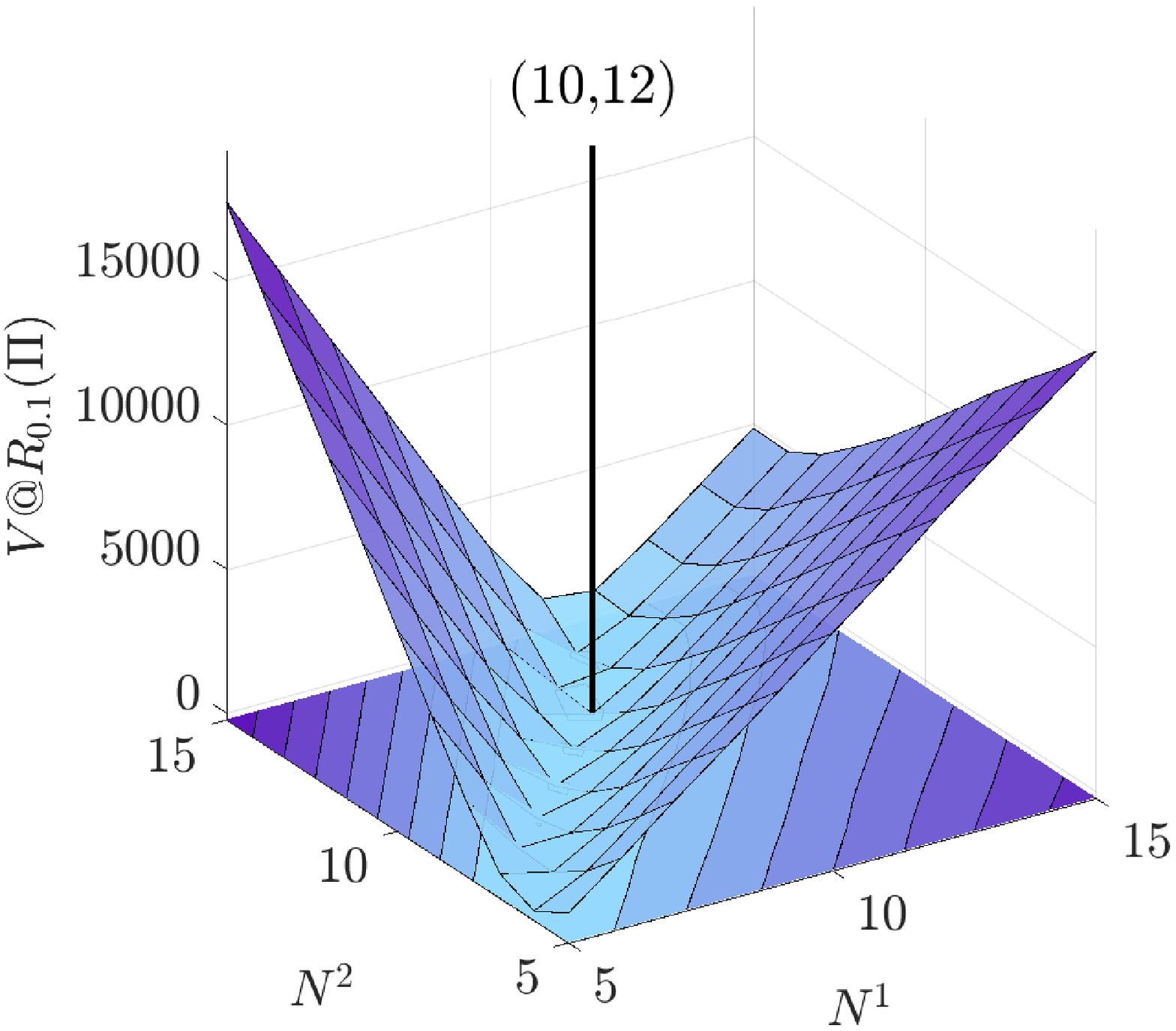}
\label{fig:Optimal_Worker_VaR01_1}
}
\subfigure[Average Value at Risk]
{
\includegraphics[width=0.46\textwidth]{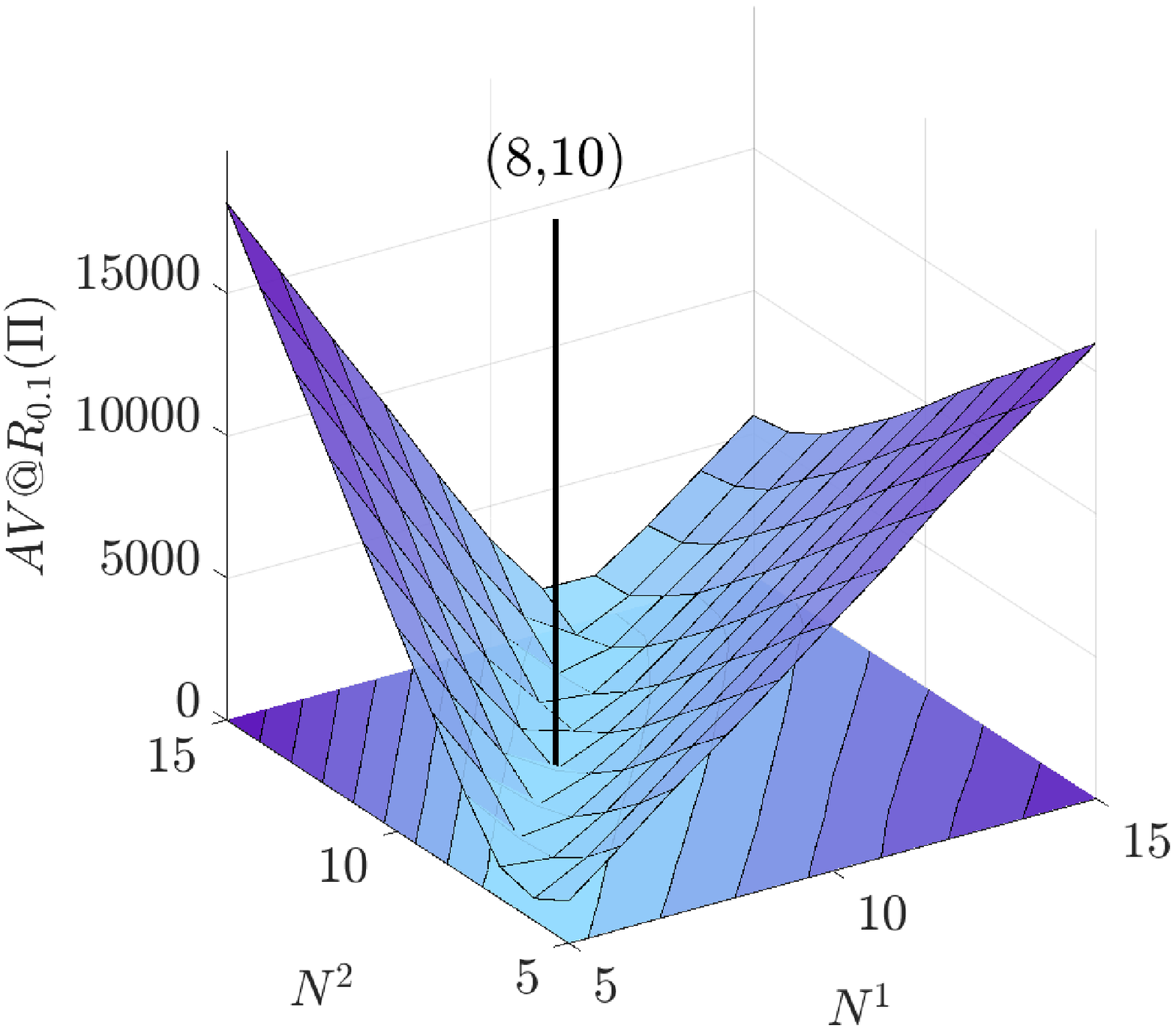}
\label{fig:Optimal_Worker_AVaR01_1}
}
\caption{Value and Average Value at Risk for different cluster sizes and level $\gl = 0.1$}
\label{fig:Profit_VaR01_1}
\end{figure}

For a level of $\gl = 0.1$, figure \ref{fig:Profit_VaR01_1} contains the Value and Average Value at Risk for this level. The cluster sizes $N^1 = 10$ and $N^2 = 12$ lead to a Value at Risk of $-0.2283 \cdot 10^3$. This can be interpreted as follows: even if we have debts of $-0.2283 \cdot 10^3$, the probability to face bankruptcy at $T = 365$ is lower than $0.1$.
The more pessimistic risk measure, Average Value at Risk, leads to a best choice $N^1 = 8$ and $N^2 = 10$ with a $\operatorname{AV@R}_{0.1}(\Pi(T))$ of $0.6255 \cdot 10^3$, i.e., we should have a surplus of $0.6255 \cdot 10^3$ for these cluster sizes.

\begin{figure}[h]
\subfigure[Value at Risk]
{
\includegraphics[width=0.46\textwidth]{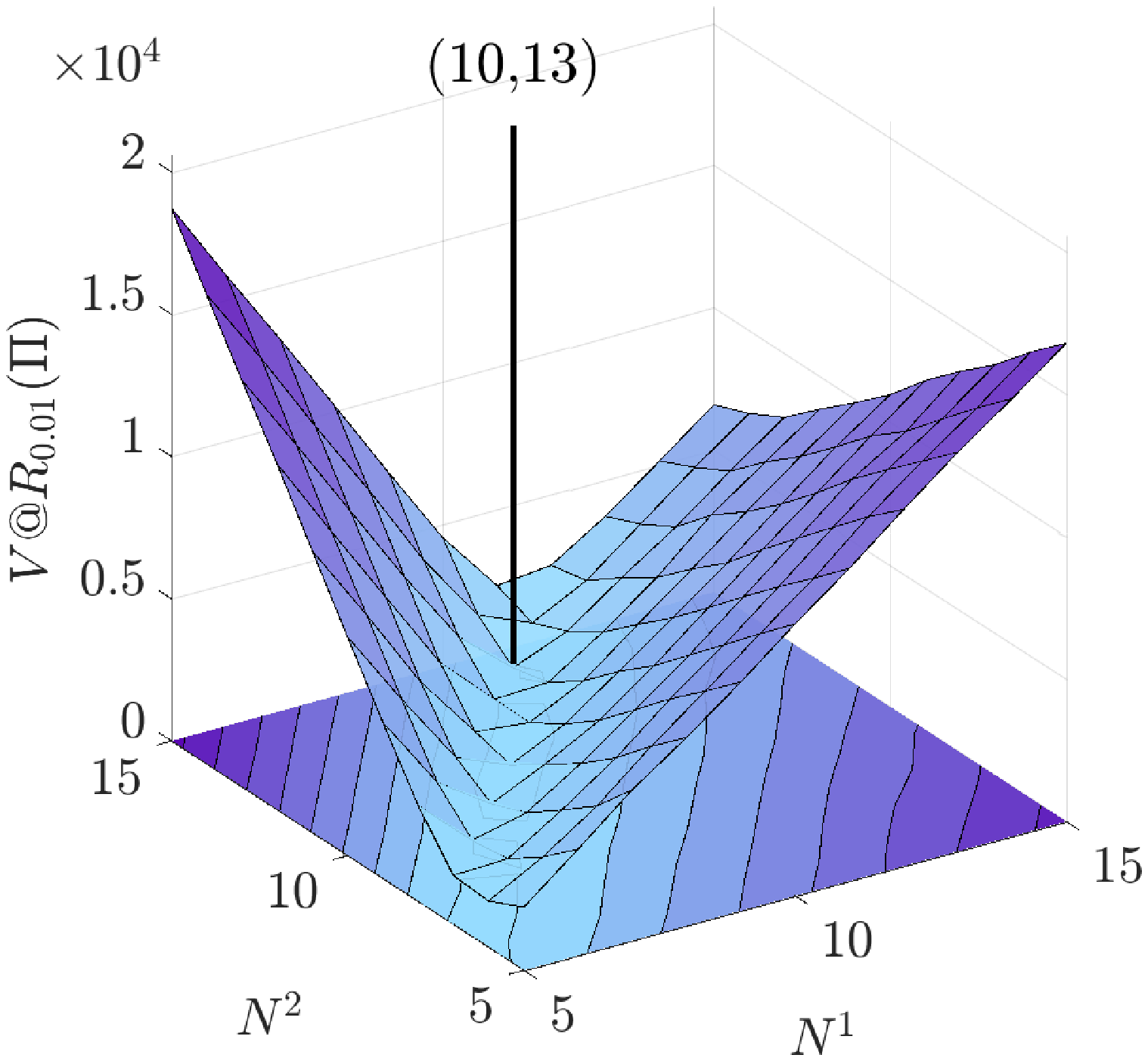}
\label{fig:Optimal_Worker_VaR001_1}
}
\subfigure[Average Value at Risk]
{
\includegraphics[width=0.46\textwidth]{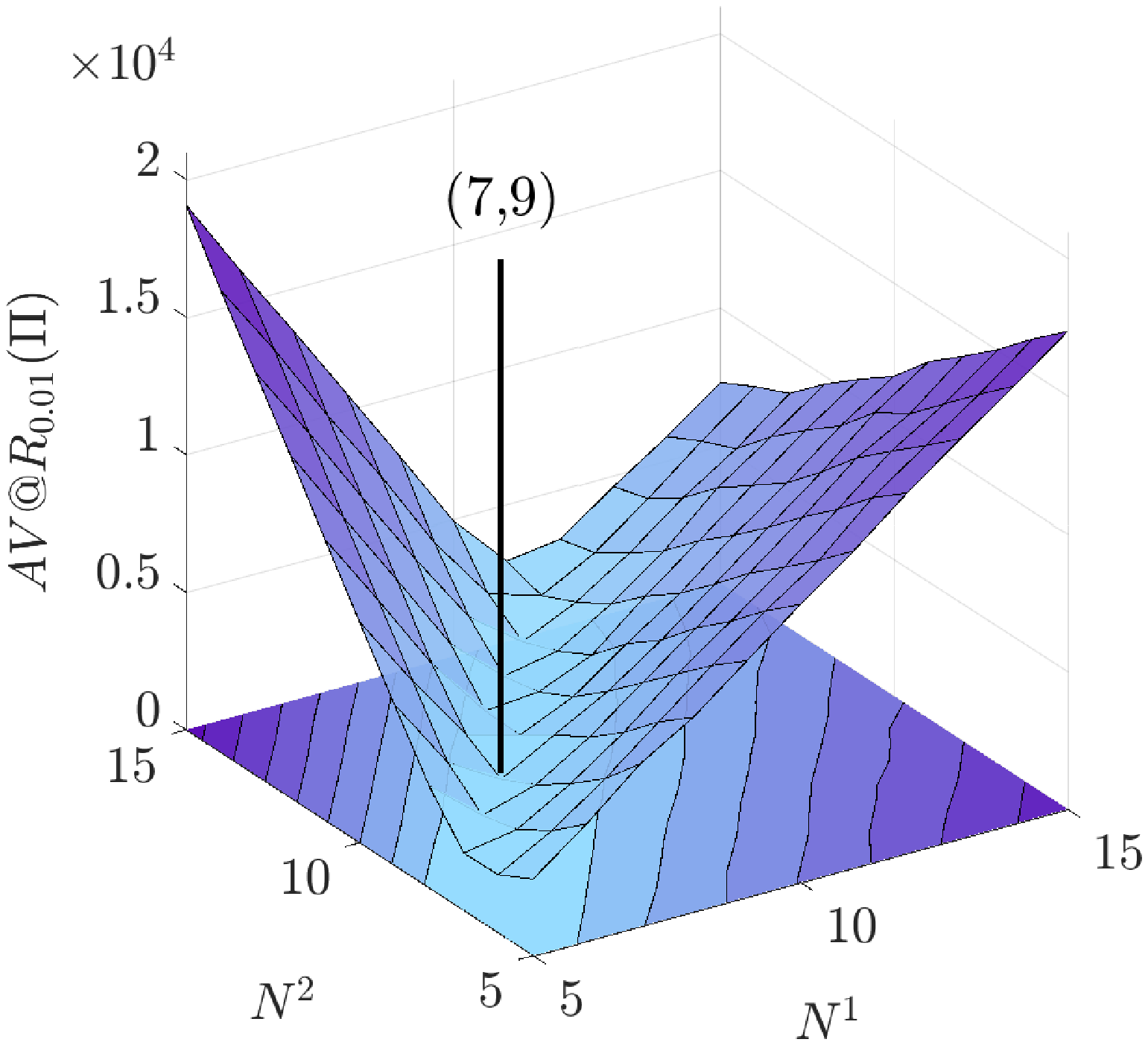}
\label{fig:Optimal_Worker_AVaR001_1}
}
\caption{Value and Average Value at Risk for different cluster sizes and level $\gl = 0.01$}
\label{fig:Profit_VaR001_1}
\end{figure}

In the case of a level of $\gl = 0.01$, the combination $N^1 = 10$ and $N^2 = 13$ yields a $\operatorname{V@R}_{0.01}(\Pi(T))$ of $1.5683 \cdot 10^3$; see figure \ref{fig:Profit_VaR001_1} (a) and table \ref{tab:optim_work}. The combination $N^1 = 7$ and $N^2 = 9$ implies an $\operatorname{AV@R}_{0.01}(\Pi(T))$ of $2.0632 \cdot 10^3$; see figure \ref{fig:Profit_VaR001_1} (b).

All introduced performance measure evaluations in table \ref{tab:optim_work} are given for reasonable cluster size choices. Obviously, the best choice for the standard deviation leads to bad performance in all the other performance measures. The allocation $N^1 = 10$ and $N^2 = 12$ implies the best solution in terms of expectation, bankruptcy probability and Value at Risk with level $0.1$. The Value and Average Value at Risk seem to work quite well here and give information about the surplus one has to hold to capture bad events and avoid a bankruptcy with a high probability (given by $1-\gl$).
One has to be careful with the interpretation here. The profit is defined as the cumulative difference of earnings and costs and does not imply detailed information at every time between zero and $T$.

From the point of optimization, the shape of the  Value at Risk and Average Value at Risk are the best, since starting from any combination we are directed into the ``valley'' of good combinations. In the case of the bankruptcy probability we have two disadvantages, namely totally flat regions and high  jumps in the values.
\vspace*{\fill}
\vspace*{\fill}
\newpage

\section{Conclusions}
We have analyzed performance measures for a stochastic production network model. It turned out that different performance measures might lead to totally different results and the choice of the appropriate measure should be a crucial point in the development of an optimization framework. 
In our examples, the measures Value at Risk and Average Value at Risk always lead to reasonable results and
a suitable curvature in the set of feasible solutions. 

Future work will deal with the formulation of rigorous optimization problems for the stochastic production network model with respect to relevant performance measures. We will study the quality of solutions and the speed of convergence to an optimum.

\section*{Acknowledgments}
This work was financially supported by BMBF project ENets (05M18VMA) and the DFG grant No. GO 1920/7-1.

\bibliographystyle{siam}
\bibliography{references}
\appendix
\section{Proof of lemma \ref{lem:sum_CTMC}}\label{app:proof1}
\begin{proof}
The proof consists of two steps. First, we prove the Markov property of $(X(t),t\geq0)$, and we compute the $\cQ$-matrix in the second part. Clearly, the stochastic process $(X(t),t\geq0)$ takes only values in the finite, discrete space $E^X:=\{0,\dots,N\}$, which together with the $\sigma$-algebra $\cE^X = \cP(E^X)$ builds a measurable space. We prove the Markov property with Dynkin's criterion \cite[theorem 10.13]{Dynkin1965} and follow the ideas of \cite[problem 4.14]{Baldi_Martingales_MC}.
We have $N$ CTMCs with state space $\{0,1\}$ given, and we construct the $N$-dimensional Markov process $(Y(t),t\geq 0)$ with $Y(t) = (X_1(t),\dots,X_N(t))$ on the product space $(\gO^N,\cA^N,P^N)$. This process takes values in the measurable space $(E^Y,\cE^Y)=(\{0,1\}^N,\cP(\{0,1\}^N))$, and we obtain the process $(X(t),t\geq 0)$ defined by the measurable surjective mapping 
\begin{align*}
\psi \colon E^Y\to E^X\text{,} \quad x \mapsto \sum_{i=1}^N x_i
\end{align*}
as $X(t) = \psi(Y(t))$. Let $(U_t,t\geq0)$ be the semigroup of Markovian kernels given by $(Y(t),t\geq 0)$, i.e., \[U_t(x,B) = P^N(Y(t+s)\in B|Y(s)=x)\] for every $x \in E^Y$, $s,t\geq 0$ and $B \in \cE^Y$. If we can show that for every $x,\tilde x \in E^Y$ with $\psi(x) = \psi(\tilde x)$ the equality \[U_t(x,\psi^{-1}(\Gamma)) = U_t(\tilde{x},\psi^{-1}(\Gamma))\]
for every $t\geq 0$ and $\Gamma \in \cE^X$ holds (Dynkin's criterion), then we know that the stochastic process $(X(t),t\geq 0)$ is a Markov process. From $x\in \psi^{-1}(\{m\})$, it follows that $x_{\sigma} \in \psi^{-1}(\{m\})$ for every $m \in E^X$ and every permutation $\sigma$ because of the structure of $\psi$. Therefore, we can compute 
\begin{align*}
&\;U_t(x,\psi^{-1}(\Gamma))\\[1ex]
=&\, \sum_{z \in \psi^{-1}(\Gamma)}P^N(Y(t) = z|Y(0) = x)\\
=&\, \sum_{z \in \psi^{-1}(\Gamma)}\prod_{i=1}^N P(X_i(t) = z_i|X_i(0) = x_i)\\
=&\, \sum_{z \in \psi^{-1}(\Gamma)}\prod_{i=1}^N P(X_1(t) = z_i|X_1(0) = x_i)\\
=&\, \sum_{z \in \psi^{-1}(\Gamma)}\prod_{i=1}^N P(X_1(t) = z_{\sigma(i)}|X_1(0) = x_{\sigma(i)})\\
=&\, \sum_{z_{\sigma} \in \psi^{-1}(\Gamma)}\prod_{i=1}^N P(X_i(t) = z_{\sigma(i)}|X_i(0) = y_i)\\[1ex]
=&\, \sum_{z_{\sigma} \in \psi^{-1}(\Gamma)}P^N(Y(t) = z_{\sigma}|Y(0) = y)\\[1ex]
=&\, U_t(y,\psi^{-1}(\Gamma))
\end{align*}
for a permutation such that $x_\sigma = y$ holds. Hence, we have the Markov property of $(X(t),t\geq 0)$ by Dynkin's criterion, and it is a CTMC.
We calculate the transition probabilities of the stochastic process $(X(t),t\geq 0)$ to obtain the $\cQ$-matrix. First, we know for every $i = 1,\dots,N$ and for every $l,m \in \{0,1\}$ that the transition probability fulfills 
\begin{align}
P(X_i(t+\Delta t) = m|X_i(t) = l) =\; \Ind_l(m)(1-\Delta t \gl_l)+\Ind_{1-l}(m) \Delta t \gl_l\;+ \lano(\Delta t) \label{eq:trans_prob1}
\end{align} 
in the limit $\Delta t\to 0$; see \cite[pp.~28]{Gross_Queueing}. 

We choose $k,j \in \{0,\dots,N\}$, and we set $C = E^N$ in the following. Then, we can write
\begin{align*}
&P(X(t+\Delta t) = k,X(t) = j)\\[1ex]
=& \sum_{\substack{\ga \in C\\|\ga| = k}}\sum_{\substack{\gb \in C\\|\gb| = j}}P\Big(X_1(t+\Delta t) = \ga_1,\dots,X_N(t+\Delta t)=\ga_N, X_1(t) = \gb_1,\dots,X_N(t) = \gb_N\Big)\\
=& \sum_{\substack{\ga \in C\\|\ga| = k}}\sum_{\substack{\gb \in C\\|\gb| = j}} \prod_{i=1}^N P(X_i(t+\Delta t) = \ga_i|X_i(t) = \gb_i)P(X_i(t) = \gb_i)
\end{align*}
using the independence of the stochastic processes and the multi-index notation \[|\ga| = \sum_{i=1}^N \ga_i.\]
At this point, we can use the transition probability \eqref{eq:trans_prob1} and merge all terms of order $\lano(\Delta t)$. This yields
\begin{align}
&P(X(t+\Delta t) = k,X(t) = j)\notag\\
&= \sum_{\substack{\ga \in C\\|\ga| = k}}\sum_{\substack{\gb \in C\\|\gb| = j}} \left(\prod_{i=1}^N \Ind_{\gb_i}(\ga_i)(1-\Delta t \gl_{\gb_i})P(X_i(t) = \gb_i)\right)\label{eq:CapDropSum1}\\
&+ \sum_{\substack{\ga \in C\\|\ga| = k}}\sum_{\substack{\gb \in C\\|\gb| = j}}\Bigg(\sum_{l=1}^N \Ind_{1-\gb_l}(\ga_l) \Delta t \gl_{\gb_l}P(X_l(t) = \gb_l) \prod_{\substack{i=1\\i\neq l}}^N \Ind_{\gb_i}(\ga_i)(1-\Delta t \gl_{\gb_i})P(X_i(t) = \gb_i)\Bigg) \label{eq:CapDropSum2}\\ &+ \lano(\Delta t)\text{.}\notag
\end{align}
To simplify the computation, we analyze both summands in the last equation separately, and we start with the first one \eqref{eq:CapDropSum1}. The product of the indicator function implies $\ga = \gb$, and we put higher orders of $\Delta t$ into $\lano(\Delta t)$ again. By doing this and observing that $k$ has to equal $j$, we see
\begin{align*}
&\sum_{\substack{\ga \in C\\|\ga| = k}}\sum_{\substack{\gb \in C\\|\gb| = j}}\left( \prod_{i=1}^N \Ind_{\gb_i}(\ga_i)(1-\Delta t \gl_{\gb_i})P(X_i(t) = \gb_i)\right)\\[1ex]
=& \,\Ind_j(k)\left(\sum_{\substack{\gb \in C\\|\gb| = j}} P(X_i(t) = \gb_i)\left(1-\Delta t \sum_{l=1}^N \gl_{\gb_l}\right)\right) + \lano(\Delta t)\\[1ex]
=& \,\Ind_j(k)\big(1-\Delta t(j\gl_1+(N-j)\gl_0)\big)P(X(t) = j)+\lano(\Delta t)\text{.}
\end{align*}
The last equality follows from
\[ \sum_{l=1}^N \gl_{\gb_l} = (j\gl_1+(N-j)\gl_0)\] since $j$ entries of $\gb$ are one and all the others are zero. Now, we analyze the second summand \eqref{eq:CapDropSum2}, which we simplify by merging terms to $\lano(\Delta t)$. We distinguish two cases, i.e., the case $k = j+1$ and the case $k = j-1$, since all remaining cases are impossible. In detail, if, e.g.,~$k = j+2$, then $\ga$ has two entries more with value one, which implies that $\ga$ and $\gb$ are different in at least two entries. Hence, the product in the second summand is zero. 
We have
\begin{align*}
&\sum_{\substack{\ga \in C\\|\ga| = k}}\sum_{\substack{\gb \in C\\|\gb| = j}}\sum_{l=1}^N \Bigg(\Ind_{1-\gb_l}(\ga_l) \Delta t \gl_{\gb_l}P(X_l(t) = \gb_l) \prod_{\substack{i=1\\i\neq l}}^N \Ind_{\gb_i}(\ga_i)(1-\Delta t \gl_{\gb_i})P(X_i(t) = \gb_i)\Bigg)\\
=&\,\Delta t\sum_{l=1}^N \sum_{\substack{\ga \in C\\|\ga| = k}}\sum_{\substack{\gb \in C\\|\gb| = j}} \Bigg(\Ind_{1-\gb_l}(\ga_l)  \gl_{\gb_l}P(X_l(t) = \gb_l) \prod_{\substack{i=1\\i\neq l}}^N \Ind_{\gb_i}(\ga_i)P(X_i(t) = \gb_i)\Bigg)+ \lano(\Delta t)\\
=&\,\Ind_{j+1}(k) \Delta t \gl_0 \sum_{l=1}^N\sum_{\substack{\gb \in C\\|\gb| = j\\\gb_l = 0}}P(X_1(t)=1)^jP(X_1(t) = 0)^{N-j}\\
&\,+\Ind_{j-1}(k) \Delta t \gl_1 \sum_{l=1}^N\sum_{\substack{\gb \in C\\|\gb| = j\\\gb_l = 1}}P(X_1(t)=1)^jP(X_1(t) = 0)^{N-j} +\lano(\Delta t)
\end{align*}
by using  that $(X_i(t),i = 1,\dots,N)$ are iid. We easily count the remaining sums with combinatorics, i.e.,~
\begin{align*}
\sum_{l=1}^N\sum_{\substack{\gb \in C\\|\gb|  = j\\\gb_l = 0}} 1 &= N\cdot \binom{N-1}{j} = (N-j)\cdot \binom{N}{j},\\
\sum_{l=1}^N\sum_{\substack{\gb \in C\\|\gb| = j\\\gb_l = 1}} 1 &= N\cdot \binom{N-1}{j-1} = j\cdot \binom{N}{j}
\end{align*}
and observe that
\begin{align*}
P(X(t) = j)= \binom{N}{j}P(X_1(t) = 1)^jP(X_1(t) = 0)^{N-j}\text{.} 
\end{align*}
Summarizing all computations yields the transition probability 
\begin{align*}
P(X(t+\Delta t) = k|X(t) = j) =&\; \Ind_j(k)(1-\Delta t(j\gl_1+(N-j)\gl_0))\\[1ex]
\;&+\Ind_{j+1}(k)\Delta t(N-j)\gl_0 \\[1ex]
\;&+\Ind_{j-1}(k)\Delta t j \gl_1 +\lano(\Delta t)
\end{align*}
as $\Delta t \to 0$ and we conclude the generator of $(X(t),t\geq 0)$ as the matrix $\cQ$ from the statement of this lemma. 
\end{proof}

\end{document}